# COINCIDENCE OF LYAPUNOV EXPONENTS FOR RANDOM WALKS IN WEAK RANDOM POTENTIALS[1]

By Markus Flury

*University of Zürich*

We investigate the free energy of nearest-neighbor random walks on $\mathbb{Z}^d$, endowed with a drift along the first axis and evolving in a nonnegative random potential given by i.i.d. random variables. Our main result concerns the ballistic regime in dimensions $d \geq 4$, at which we show that quenched and annealed Lyapunov exponents are equal as soon as the strength of the potential is small enough.

**1. Introduction and main results.**

1.1. *Random walk in random potential.* Let $\mathcal{S} = (S(n))_{n \in \mathbb{N}_0}$ be a nearest-neighbor random walk on the lattice $\mathbb{Z}^d$, with start at the origin and drift $h$ in the direction of the first axis. We suppose $\mathcal{S}$ to be defined on a probability space $(\Omega, \mathcal{F}, P_h)$ and we denote by $E_h$ the associated expectation. Such a random process is characterized by the distributions of its finite-step subpaths

$$\mathcal{S}[n] \stackrel{\text{def}}{=} (S(0), \ldots, S(n)), \qquad n \in \mathbb{N}.$$

In the nondrifting case $h = 0$, these distributions are uniform on the nearest-neighbor paths in $\mathbb{Z}^d$. That is, for $n \in \mathbb{N}$ and $x_0$ the origin, we have

$$P_0[\mathcal{S}[n] = (x_0, \ldots, x_n)] = \left(\frac{1}{2d}\right)^n$$

for all $x_1, \ldots, x_n \in \mathbb{Z}^d$ such that $\|x_i - x_{i-1}\| = 1$ for $i = 1, \ldots, n$, the probability being zero elsewhere. The case of a nonvanishing drift is related back

Received October 2006; revised September 2007.
[1]Supported by the Swiss National Fund under contract nos. 20-55648.98, 20-63798.00 and 20-100536/1.
*AMS 2000 subject classifications.* Primary 60K37; secondary 34D08, 60K35.
*Key words and phrases.* Random walk, random potential, Lyapunov exponents, interacting path potential.

to the nondrifting case by means of the density function

$$(1.1) \qquad \frac{dP_h \mathcal{S}[n]^{-1}}{dP_0 \mathcal{S}[n]^{-1}} = \frac{\exp(h \cdot S_1(n))}{E_0[\exp(h \cdot S_1(n))]}, \qquad n \in \mathbb{N},$$

with $S_1$ denoting the first component of $S$.

The random walk $\mathcal{S}$ is a *Markov chain* with *independent growths*, which means the following. Suppose $m, n \in \mathbb{N}_0$ and set

$$\mathcal{S}[m, n] \stackrel{\text{def}}{=} (S(m), \ldots, S(n)).$$

When $x_0, \ldots, x_n \in \mathbb{Z}^d$ fulfill $P_h[\mathcal{S}[n] = (x_0, \ldots, x_m)] > 0$, we have

$$P_h[\mathcal{S}[m, n] = (x_m, \ldots, x_n) \mid \mathcal{S}[m] = (x_0, \ldots, x_m)]$$
$$= P_h[\mathcal{S}[n - m] = (x_m - x_m, \ldots, x_n - x_m)].$$

We will constantly make use of this property and will simply refer to it as the *Markov property* (committing a slight abuse of standard terminology).

In addition to the influence of the drift, we want $\mathcal{S}$ to underlie the influence of a random potential on the lattice. To this end, let $\mathbb{V} = (V_x)_{x \in \mathbb{Z}^d}$ be a family of independent, identically distributed random variables, independent of the random walk itself, with $\operatorname{ess\,inf} V_x = 0$ and $\mathbb{E} V_x^d < \infty$. To avoid trivialities, we also assume $\mathbb{P}[V_x = 0] < 1$.

Using the random potential $\mathbb{V}$, we are now able to introduce path measures for the random walk. Thereby, we distinguish between the so-called *quenched* setting, where the path measure depends on the concrete realization of the potential $\mathbb{V}$, and the *annealed* setting, where the measure depends on averaged values of the potential only.

The *quenched path potential* is given by

$$\Phi^{\text{qu}}_{\mathbb{V},\beta}(N) \stackrel{\text{def}}{=} \beta \sum_{n=1}^{N} V_{S(n)}, \qquad N \in \mathbb{N},$$

where the so-called *inverted temperature* $\beta \geq 0$ is a parameter for the strength of the potential. The *quenched path measure* is defined by means of the density function

$$\frac{dQ^{\text{qu}}_{\mathbb{V},h,\beta,N}}{dP_h} \stackrel{\text{def}}{=} \frac{\exp(-\Phi^{\text{qu}}_{\mathbb{V},\beta}(N))}{Z^{\text{qu}}_{\mathbb{V},h,\beta,N}}, \qquad N \in \mathbb{N},$$

where the normalization

$$Z^{\text{qu}}_{\mathbb{V},h,\beta,N} \stackrel{\text{def}}{=} E_h[\exp(-\Phi^{\text{qu}}_{\mathbb{V},\beta}(N))], \qquad N \in \mathbb{N},$$

is called the *quenched partition function*. The quenched setting defines a discrete-time model for a particle moving in a random medium. Here, the path measure is itself random, the randomness coming from the random environment $\mathbb{V}$. Under a concrete realization of the path measure, the walker jumps from site to site, thereby trying to stay in regions where the potential takes on small values. The drift, however, implies a restriction in the search for such an "optimal strategy" by imposing a particular direction on the walk.

The *annealed path measure* is defined by means of the density function

$$\frac{dQ^{\mathrm{an}}_{h,\beta,N}}{dP_h} \stackrel{\mathrm{def}}{=} \frac{\mathbb{E}\exp(-\Phi^{\mathrm{qu}}_{\mathbb{V},\beta}(N))}{\mathbb{E}Z^{\mathrm{qu}}_{\mathbb{V},h,\beta,N}}, \qquad N \in \mathbb{N},$$

and $\mathbb{E}Z^{\mathrm{qu}}_{\mathbb{V},h,\beta,N}$ is called the *annealed partition function*. While our main interest lies in the quenched setting, the annealed model no longer depends on the realizations of the environment and is thus easier to handle. A walker under the annealed measure finds himself in a similar situation as in the quenched setting. To see this, observe that the quenched potential can be expressed by

$$\Phi^{\mathrm{qu}}_{\mathbb{V},\beta}(N) = \beta \sum_{x \in \mathbb{Z}^d} \ell_x(N) V_x,$$

where

$$\ell_x(N) \stackrel{\mathrm{def}}{=} \sum_{n=1}^{N} 1_{\{S(n)=x\}}$$

denotes the number of visits to the site $x \in \mathbb{Z}^d$ by the $N$-step random walk $\mathcal{S}[1,N]$. An *annealed path potential* is given by

$$\Phi^{\mathrm{an}}_{\beta}(N) \stackrel{\mathrm{def}}{=} \sum_{x \in \mathbb{Z}^d} \varphi^{\mathrm{an}}_{\beta}(\ell_x(N)), \qquad N \in \mathbb{N},$$

where

(1.2) $$\varphi^{\mathrm{an}}_{\beta}(t) \stackrel{\mathrm{def}}{=} -\log \mathbb{E}\exp(-t\beta V_x), \qquad t \in \mathbb{R}^+,$$

is a nonnegative function which is concave increasing by the Hölder inequality. Now, by the independence assumption on $\mathbb{V}$, it is easily seen that

$$\frac{dQ^{\mathrm{an}}_{h,\beta,N}}{dP_h} = \frac{\exp(-\Phi^{\mathrm{an}}_{\beta}(N))}{Z^{\mathrm{an}}_{h,\beta,N}},$$

where the normalizing constant

$$Z^{\mathrm{an}}_{h,\beta,N} \stackrel{\mathrm{def}}{=} E_h[\exp(-\Phi^{\mathrm{an}}_{\beta}(N))]$$

equals the annealed partition function $\mathbb{E}Z^{\mathrm{qu}}_{\mathbb{V},h,\beta,N}$. By the concavity of $\varphi^{\mathrm{an}}_{\beta}$, the more often the random walk intersects its own path, the smaller the potential $\Phi^{\mathrm{an}}_{\beta}$. Therefore, on the one hand, it is convenient for the walker to return to places he has already visited, while on the other hand, he is urged to proceed in the direction of the drift.

In a similar model in a continuous setting, namely Brownian motion in a Poissonian potential, the contrary influence of drift and potential on the long-time behavior of the walk was first studied by A. S. Sznitman. By means of the powerful method of enlargement of obstacles, he established a precise picture in both quenched and annealed settings (see Chapter 5 of his book [12]). Among his results there is an accurate description of a phase transition from *localization* for large $\beta$ (or small $h$) to *delocalization* for small $\beta$ (or large $h$). In the delocalized phase, the random walk is *ballistic*, that is, the displacement of $S(N)$ from the origin grows of order $O(N)$, while in the localized phase, the walk behaves *sub-ballistic*, that is, the displacement is of order $o(N)$. The analogous results for the discrete setting have been established by Zerner in [15] and Flury in [7].

1.2. *Lyapunov exponents.* The above results on the transition from sub-ballistic to ballistic behavior are based on large deviation principles for the random walk under the path measures and on phase transitions for the *quenched* and *annealed free energies*

$$\log Z^{\mathrm{qu}}_{\mathbb{V},h,\beta,N} \quad \text{and} \quad \log Z^{\mathrm{an}}_{h,\beta,N}.$$

The free energies are important values for the study of the path measures. The moments of the path potential under these measures, for instance, may be evaluated by differentiating the free energies with respect to the inverted temperature. For a more direct motivation in the context of random branching processes, and a thorough study of the one-dimensional case, we refer to [8] by Greven and den Hollander.

The main subject of the present article is the long-time behavior of the free energies. We first deal with the phase transitions in this behavior, as established in [7]. The associated phase diagrams coincide with the ones for the random walk itself. They are characterized by values from the so-called *point-to-hyperplane* setting. For $h > 0$, $\beta \geq 0$ and $L \in \mathbb{N}$, we set

$$\overline{Z}^{\mathrm{qu}}_{\mathbb{V},h,\beta,L} \stackrel{\mathrm{def}}{=} \sum_{n=1}^{\infty} E_h[\exp(-\Phi^{\mathrm{qu}}_{\mathbb{V},\beta}(n)); \{S_1(n) = L\}],$$

$$\overline{Z}^{\mathrm{an}}_{h,\beta,L} \stackrel{\mathrm{def}}{=} \sum_{n=1}^{\infty} E_h[\exp(-\Phi^{\mathrm{an}}_{\beta}(n)); \{S_1(n) = L\}].$$

THEOREM A. *For any $\beta \geq 0$, there are continuous, nonnegative functions $\overline{m}^{\mathrm{qu}}(\cdot, \beta)$ and $\overline{m}^{\mathrm{an}}(\cdot, \beta)$ on $\mathbb{R}^+$ such that*

$$\overline{m}^{\mathrm{qu}}(h, \beta) = -\lim_{L \to \infty} \frac{1}{L} \log \overline{Z}^{\mathrm{qu}}_{\mathbb{V},h,\beta,L},$$

$$\overline{m}^{\mathrm{an}}(h,\beta) = -\lim_{L\to\infty} \frac{1}{L} \log \overline{Z}^{\mathrm{an}}_{h,\beta,L}$$

for all $h > 0$, as well as continuous, nonnegative functions $m^{\mathrm{qu}}(\cdot,\beta)$ and $m^{\mathrm{an}}(\cdot,\beta)$ on $\mathbb{R}^+$, the so-called quenched and annealed Lyapunov exponents, such that

$$m^{\mathrm{qu}}(h,\beta) = -\lim_{N\to\infty} \frac{1}{N} \log Z^{\mathrm{qu}}_{\mathbb{V},h,\beta,N},$$

$$m^{\mathrm{an}}(h,\beta) = -\lim_{N\to\infty} \frac{1}{N} \log Z^{\mathrm{an}}_{h,\beta,N}$$

for all $h \geq 0$, where the convergence in the quenched setting is $\mathbb{P}$-almost surely and in $\mathcal{L}_1(\mathbb{P})$, and where the limits no longer depend on the realizations of $\mathbb{V}$. Moreover, for all $h, \beta \geq 0$, we have

$$m^{\mathrm{qu}}(h,\beta) = \begin{cases} \lambda_h, & \text{if } \overline{m}^{\mathrm{qu}}(0,\beta) \geq h, \\ \lambda_h - \lambda_{\bar{h}^{\mathrm{qu}}(h,\beta)}, & \text{if } \overline{m}^{\mathrm{qu}}(0,\beta) < h, \end{cases}$$

$$m^{\mathrm{an}}(h,\beta) = \begin{cases} \lambda_h, & \text{if } \overline{m}^{\mathrm{an}}(0,\beta) \geq h, \\ \lambda_h - \lambda_{\bar{h}^{\mathrm{an}}(h,\beta)}, & \text{if } \overline{m}^{\mathrm{an}}(0,\beta) < h, \end{cases}$$

where $\bar{h}^{\mathrm{qu}}(h,\beta) > 0$ and $\bar{h}^{\mathrm{an}}(h,\beta) > 0$ are determined by

$$\overline{m}^{\mathrm{qu}}(\bar{h}^{\mathrm{qu}}(h,\beta),\beta) = h - \bar{h}^{\mathrm{qu}}(h,\beta),$$

$$\overline{m}^{\mathrm{an}}(\bar{h}^{\mathrm{an}}(h,\beta),\beta) = h - \bar{h}^{\mathrm{an}}(h,\beta)$$

and where $\lambda_h \stackrel{\mathrm{def}}{=} \log E_0[\exp(h \cdot S_1(1))]$.

Theorem A is proved in [7] for drifts in arbitrary directions and for a more general annealed potential which we introduce at the beginning of Section 2. With regard to the difference in notation, observe that, with $e_1$ being the first unit vector,

$$m^{\mathrm{qu}}(h,\beta) = \lambda_h - \lim_{n\to\infty} \frac{1}{n} \log Z^{h \cdot e_1}_{n,\omega}, \qquad \mathbb{P}\text{-a.s.},$$

$$m^{\mathrm{an}}(h,\beta) = \lambda_h - \lim_{n\to\infty} \frac{1}{n} \log Z^{h \cdot e_1}_n$$

and that by Corollary C of [7],

$$\overline{m}^{\mathrm{qu}}(\bar{h},\beta) = \frac{1}{\alpha^*_{\lambda_{\bar{h} \cdot e_1}}(e_1)} - \bar{h},$$

$$\overline{m}^{\mathrm{an}}(\bar{h},\beta) = \frac{1}{\beta^*_{\lambda_{\bar{h} \cdot e_1}}(e_1)} - \bar{h},$$

where the notation on the right-hand sides is from [7] (with potential $\mathbb{V}_\beta = \{\beta V_x\}$ and $\varphi = \varphi_{\mathbb{V}_\beta}$). For the latter equalities, observe also that the fact that the random walk is already stopped at its first entrance into the hyperplane has no effect on the limits in Corollary C of [7] (as will become clear in Section 2.1).

In accordance with the long-time behavior of the random walk itself, in Theorem A, we have the following picture for the behavior of the free energies: in the sub-ballistic regime, the walker remains near the origin in the annealed case and in regions with small potential in the quenched case. Therefore, since

$$\lim_{t \to \infty} \frac{\varphi_\beta^{\mathrm{an}}(t)}{t} = \operatorname{ess\,inf} V_x = 0,$$

the contribution from the potential then vanishes when $N$ becomes large. What remains is the probability of staying in an only slow-growing region, contributing the value

$$\lambda_h = \lim_{n \to \infty} \frac{1}{n} \log E_0[\exp(h \cdot S_1(n))].$$

In the ballistic regime, on the other hand, the walk obeys the drift and dislocates with a nonvanishing velocity. As a consequence, the path potential and the "spatial part" of the density for the drift must not be neglected, as they contribute the subtraction term $\lambda_{\bar{h}^{\mathrm{qu}}(h,\beta)}$, respectively $\lambda_{\bar{h}^{\mathrm{an}}(h,\beta)}$, to the corresponding Lyapunov exponent (see [7] for a rigorous interpretation of this last point).

1.3. *Main results and preliminaries.* Our first new result concerns the simpler annealed setting. For $h > 0$, the critical parameter $\beta_{\mathrm{c}}^{\mathrm{an}}(h)$ for the phase transition is given by

$$\overline{m}^{\mathrm{an}}(0, \beta_{\mathrm{c}}^{\mathrm{an}}(h)) = h,$$

where existence and uniqueness of $\beta_{\mathrm{c}}^{\mathrm{an}}(h)$ will be explained in Remark 2.10 of Section 2.1.

THEOREM B. (a) *For any $h, \beta \geq 0$, we have*

$$Z_{h,\beta,N}^{\mathrm{an}} \leq \exp(-m^{\mathrm{an}}(h,\beta)N), \qquad N \in \mathbb{N},$$

*and $m^{\mathrm{an}}$ is continuous on $\mathbb{R}^+ \times \mathbb{R}^+$. Moreover, for any $h > 0$ and $\beta_0 < \beta_{\mathrm{c}}^{\mathrm{an}}(h)$, there exists $K_{h,\beta_0} < \infty$ such that for $\beta \leq \beta_0$, we have*

$$Z_{h,\beta,N}^{\mathrm{an}} \geq \frac{\exp(-m^{\mathrm{an}}(h,\beta)N)}{K_{h,\beta_0}}, \qquad N \in \mathbb{N}.$$

(b) *$m^{\mathrm{an}}$ is analytic on the open set $\{(h,\beta) \in (0,\infty)^2 : \beta \neq \beta_{\mathrm{c}}^{\mathrm{an}}(h)\}$.*

Part (a) of Theorem B is established in Section 3.1, essentially by subadditivity arguments. The sub-ballistic part in (b) is a straightforward consequence of Theorem A. The more complicated ballistic part is proved in Section 3.2 by renewal techniques.

The next theorem is the main result of this article. It concerns dimensions $d \geq 4$ and nonvanishing drifts $h > 0$. It states that quenched and annealed Lyapunov exponents coincide once the strength of the potential is chosen to be small enough.

THEOREM C. *Suppose $d \geq 4$ and $h > 0$. There then exists $\beta_0 > 0$ such that*

$$m^{\mathrm{qu}}(h,\beta) = m^{\mathrm{an}}(h,\beta)$$

*for all $\beta \leq \beta_0$. Moreover, when $V_x$ is essentially bounded, there exists $K_{\mathrm{fr.e.}} < \infty$ such that*

$$\mathbb{E}|\log Z^{\mathrm{qu}}_{\mathbb{V},h,\beta,N} - \log Z^{\mathrm{an}}_{h,\beta,N}| \leq K_{\mathrm{fr.e.}}(1 + \beta\sqrt{N})$$

*for all $N \in \mathbb{N}$ and $\beta \leq \beta_0$.*

Coincidence of Lyapunov exponents has been conjectured by Sznitman in [12]. It emerged from the fact that a similar result is true for the much simpler case of directed polymers in random potentials. There, the random walk $(S(n))_{n \in \mathbb{N}}$ is replaced by $((\xi(n),n))_{n \in \mathbb{N}}$, where $(\xi(n))_{n \in \mathbb{N}}$ is a standard $d$-dimensional walk. The coincidence of quenched and annealed Lyapunov exponents for $d \geq 3$ and small disorder was first proven by Imbrie and Spencer in [10], using cluster expansion techniques, and then by Bolthausen in [2] and Albeverio and Zhou in [1], using martingale techniques. Martingale arguments have also been used in the more recent work on directed polymers in [5] and [3].

The situation considered here is much more delicate and, unfortunately, it does not seem possible to implement martingale techniques. We therefore take recourse to different methods, mainly renewal techniques and arguments from Ornstein–Zernike theory.

The crucial result for the proof of Theorem C is an estimate on the second moment of the quenched partition function.

THEOREM D. *Suppose $d \geq 4$ and $h > 0$. There are then $\beta_0 > 0$ and $K_{\mathrm{s.m.}} < \infty$ such that*

$$\mathbb{E}(Z^{\mathrm{qu}}_{\mathbb{V},h,\beta,N})^2 \leq K_{\mathrm{s.m.}}(Z^{\mathrm{an}}_{h,\beta,N})^2$$

*for all $N \in \mathbb{N}$ and $\beta \leq \beta_0$.*

In order to achieve a heuristic understanding of Theorem D, we consider two independent copies $\mathcal{S}^1 = (S^1(n))_{n \in \mathbb{N}_0}$ and $\mathcal{S}^2 = (S^2(n))_{n \in \mathbb{N}_0}$ of the random walk $\mathcal{S}$. For $x \in \mathbb{Z}^d$, we set

$$\ell_x^1(N) \stackrel{\text{def}}{=} \sum_{n=1}^N 1_{\{S^1(n) = x\}} \quad \text{and} \quad \ell_x^2(N) \stackrel{\text{def}}{=} \sum_{n=1}^N 1_{\{S^2(n) = x\}}.$$

Recall that we have $Z_{h,\beta,N}^{\text{an}} = \mathbb{E} Z_{\mathbb{V},h,\beta,N}^{\text{qu}}$, by the independence assumption on the potential. In a similar way, and by the independence of $\mathcal{S}^1$ and $\mathcal{S}^2$, we obtain

$$\mathbb{E}(Z_{\mathbb{V},h,\beta,N}^{\text{qu}})^2 = E_h\left[\exp\left(-\sum_{x \in \mathbb{Z}^d} \varphi_\beta^{\text{an}}(\ell_x^1(N) + \ell_x^2(N))\right)\right],$$

$$(\mathbb{E} Z_{\mathbb{V},h,\beta,N}^{\text{qu}})^2 = E_h\left[\exp\left(-\sum_{x \in \mathbb{Z}^d} \varphi_\beta^{\text{an}}(\ell_x^1(N)) + \varphi_\beta^{\text{an}}(\ell_x^2(N))\right)\right],$$

where $E_h$ denotes the expectation with respect to the product measure $P_h \otimes P_h$. With the further notation

(1.3) $\quad \Psi_{\beta,N} \stackrel{\text{def}}{=} \sum_{x \in \mathbb{Z}^d} \varphi_\beta^{\text{an}}(\ell_x^1(N)) + \varphi_\beta^{\text{an}}(\ell_x^2(N)) - \varphi_\beta^{\text{an}}(\ell_x^1(N) + \ell_x^2(N))$

and with $E_{h,\beta,N}^{\text{an}}$ the expectation with respect to the annealed product path measure $Q_{h,\beta,N}^{\text{an}} \otimes Q_{h,\beta,N}^{\text{an}}$, we thus have

(1.4) $\quad\quad\quad\quad\quad \dfrac{\mathbb{E}(Z_{\mathbb{V},h,\beta,N}^{\text{qu}})^2}{(\mathbb{E} Z_{\mathbb{V},h,\beta,N}^{\text{qu}})^2} = E_{h,\beta,N}^{\text{an}}[\exp(\Psi_{\beta,N})].$

Observe further that the only nonvanishing summands in (1.3) are the ones associated to those $x \in \mathbb{Z}^d$ that are visited by both random walks up to time $N$. From the concavity of $\varphi_\beta^{\text{an}}$, we therefore obtain

(1.5) $\quad\quad\quad\quad\quad \Psi_{\beta,N} \leq \varphi_\beta^{\text{an}}(1) \sum_{x \in R^1(N) \cap R^2(N)} \ell_x^1(N) + \ell_x^2(N),$

where $R^j(N) \stackrel{\text{def}}{=} \{S^j(n) : n = 1, \ldots, N\}$ for $j = 1, 2$. This finally gives us the following picture of the situation: in the ballistic regime, the random walks $\mathcal{S}^1$ and $\mathcal{S}^2$ under the annealed path measure obey the drift and evolve in the direction of the first axis. While they do that, one expects them to move away from each other in the $(d-1)$-dimensional "vertical" direction as soon as the dimension of the lattice is large enough. The condition $d \geq 4$ appears to be the right one since the "vertical distance" is then transient. As a consequence, the paths of $\mathcal{S}^1$ and $\mathcal{S}^2$ are supposed to intersect only finitely many times. For $\beta$ small enough, the right-hand side of (1.4) should then stay bounded as $N \to \infty$ because of (1.5) and $\lim_{\beta \downarrow 0} \varphi_\beta^{\mathrm{an}}(1) = 0$.

The equality of the Lyapunov exponents, knowing that the quenched free energy is deterministic, is obtained from Theorem D by rather elementary methods. On the other hand, the estimate of the speed of convergence for the free energy has a more complicated derivation, requiring a concentration inequality for the quenched free energy. In this particular model, the usual concentration estimate is not sharp enough. For this reason, we replace $Z_{\mathbb{V},h,\beta,N}^{\mathrm{qu}}$ by a modified partition function. For $h > 0$, $\beta \geq 0$, $k \in \mathbb{R}^+$ and $N \in \mathbb{N}$, we define

$$Z_{\mathbb{V},h,\beta,k,N}^{\mathrm{qu}} \stackrel{\mathrm{def}}{=} E_h\left[e^{-\beta \sum_{x \in \mathbb{Z}^d} \ell_x(N) V_x(\omega)}; \sum_{x \in \mathbb{Z}^d} \ell_x(N)^2 \leq kN\right].$$

The justification for such a replacement is given in the following lemma, which is proved in Section 3.3.

LEMMA E. *Suppose $h > 0$ and $\beta_0 < \beta_{\mathrm{c}}^{\mathrm{an}}(h)$, and let $K_{h,\beta_0}$ be chosen according to Theorem B. For any $\varepsilon < 1/K_{h,\beta_0}$, there then exists $k_\varepsilon < \infty$ such that*

$$\mathbb{E} Z_{\mathbb{V},h,\beta,k_\varepsilon,N}^{\mathrm{qu}} \geq \varepsilon Z_{h,\beta,N}^{\mathrm{an}}$$

*for all $N \in \mathbb{N}$ and $\beta \leq \beta_0$.*

By Theorem A, Theorem D and Lemma E, we now have the means to prove the coincidence of the Lyapunov exponents and to estimate the speed of convergence for the free energy.

PROOF OF THEOREM C. By the quenched part of Theorem A, we have

$$\lim_{N \to \infty} \frac{1}{N} \mathbb{E} \log Z_{\mathbb{V},h,\beta,N}^{\mathrm{qu}} = -m^{\mathrm{qu}}(h,\beta)$$

and thus $m^{\mathrm{qu}}(h,\beta) \geq m^{\mathrm{an}}(h,\beta)$, by Jensen's inequality. In order to obtain the inverted estimate, observe that

$$\mathbb{P}\left[\lim_{N \to \infty} \frac{1}{N} \log Z_{\mathbb{V},h,\beta,N}^{\mathrm{qu}} \geq -m^{\mathrm{an}}(h,\beta)\right] \geq \liminf_{N \to \infty} \mathbb{P}\left[Z_{\mathbb{V},h,\beta,N}^{\mathrm{qu}} \geq \frac{1}{2} \mathbb{E} Z_{\mathbb{V},h,\beta,N}^{\mathrm{qu}}\right]$$

and that the left-hand side is either one or zero since the limit is deterministic (again by Theorem A). The Schwarz inequality now implies

$$\mathbb{E} Z_{\mathbb{V},h,\beta,N}^{\mathrm{qu}} \leq \tfrac{1}{2} \mathbb{E} Z_{\mathbb{V},h,\beta,N}^{\mathrm{qu}} + (\mathbb{E}(Z_{\mathbb{V},h,\beta,N}^{\mathrm{qu}})^2)^{1/2} \mathbb{P}[Z_{\mathbb{V},h,\beta,N}^{\mathrm{qu}} \geq \tfrac{1}{2} \mathbb{E} Z_{\mathbb{V},h,\beta,N}^{\mathrm{qu}}]^{1/2},$$

which leads to the Paley–Zigmund inequality

(1.6) $$\mathbb{P}\left[Z_{\mathbb{V},h,\beta,N}^{\mathrm{qu}} \geq \frac{1}{2} \mathbb{E} Z_{\mathbb{V},h,\beta,N}^{\mathrm{qu}}\right] \geq \frac{1}{4} \frac{(\mathbb{E} Z_{\mathbb{V},h,\beta,N}^{\mathrm{qu}})^2}{\mathbb{E}(Z_{\mathbb{V},h,\beta,N}^{\mathrm{qu}})^2}.$$

The lower estimate for $m^{\mathrm{qu}}(h,\beta)$ thus follows from Theorem D.

We proceed to the estimate for the speed of convergence. Assume that

$$V \stackrel{\mathrm{def}}{=} \operatorname{ess\,sup} V_x < \infty.$$

We first investigate the modified partition function $Z^{\mathrm{qu}}_{\mathbb{V},h,\beta,k,N}$. For $N \in \mathbb{N}$, let $M^{\mathrm{qu}}_{h,\beta,k,N}$ be a median of $\log Z^{\mathrm{qu}}_{\mathbb{V},h,\beta,k,N}$, that is, a real number $M^{\mathrm{qu}}_{h,\beta,k,N}$ with

$$\mathbb{P}[\log Z^{\mathrm{qu}}_{\mathbb{V},h,\beta,k,N} \leq M^{\mathrm{qu}}_{h,\beta,k,N}] \geq \tfrac{1}{2} \text{ and } \mathbb{P}[\log Z^{\mathrm{qu}}_{\mathbb{V},h,\beta,k,N} \geq M^{\mathrm{qu}}_{h,\beta,k,N}] \geq \tfrac{1}{2}.$$

Also, let $B_N \stackrel{\mathrm{def}}{=} \{x \in \mathbb{Z}^d : \|x\|_1 \leq N\}$ and let $f : [-1,1]^{B_N} \to \mathbb{R}$ be given by

$$f_k(v) \stackrel{\mathrm{def}}{=} \log E_h\left[e^{-\beta V \sum_{x \in B_N} \ell_x(N) v_x} ; \sum_{x \in B_N} \ell_x(N)^2 \leq kN\right]$$

for $v = (v_x)_{x \in B_N} \in [-1,1]^{B_N}$. We then obviously have

$$\log Z^{\mathrm{qu}}_{\mathbb{V},h,\beta,k,N} = f_k \circ (V_x/V)_{x \in B_N}.$$

Since the function $f$ is convex by the Hölder inequality, the sets $f^{-1}((-\infty, a])$ for $a \in \mathbb{R}$ are also convex. In addition, for any $v, w \in [-1,1]^{B_N}$, we have

$$|f_k(v) - f_k(w)|$$

$$\leq \sup\left\{\beta V \sum_{x \in B_N} \ell_x |v_x - w_x| : \ell \in \mathbb{N}_0^{B_N} \text{ with } \sum_{x \in B_N} \ell_x^2 \leq kN\right\}$$

$$\leq \beta V \sqrt{kN} \left(\sum_{x \in B_N} |v_x - w_x|^2\right)^{1/2}$$

by the Cauchy–Schwarz inequality for sums. This means that $f_k$ is Lipschitz continuous with Lipschitz constant at most $\beta V \sqrt{kN}$. We can thus apply Theorem 6.6 of [13] to obtain the concentration inequality

$$(1.7) \quad \mathbb{P}[|\log Z^{\mathrm{qu}}_{\mathbb{V},h,\beta,k,N} - M^{\mathrm{qu}}_{h,\beta,k,N}| \geq t] \leq 4\exp\left(\frac{-t^2}{16\beta^2 V^2 kN}\right), \quad t \in \mathbb{R}^+.$$

Next, we will find an estimate for

$$\mathbb{E}|\log Z^{\mathrm{qu}}_{\mathbb{V},h,\beta,k,N} - \log \mathbb{E}Z^{\mathrm{qu}}_{\mathbb{V},h,\beta,k,N}|$$

by adding and subtracting $M^{\mathrm{qu}}_{h,\beta,k,N}$. To this end, observe that (1.7) implies that

$$\mathbb{E}|\log Z^{\mathrm{qu}}_{\mathbb{V},h,\beta,k,N} - M^{\mathrm{qu}}_{h,\beta,k,N}| = \int_0^\infty \mathbb{P}[|\log Z^{\mathrm{qu}}_{\mathbb{V},h,\beta,k,N} - M^{\mathrm{qu}}_{h,\beta,k,N}| > t]\,dt$$

$$(1.8) \qquad\qquad \leq 4\int_0^\infty \exp\left(\frac{-t^2}{16\beta^2 kV^2 N}\right) dt$$

$$= 8V\beta\sqrt{k\pi N}.$$

It remains to find an estimate for $|\log \mathbb{E} Z^{\mathrm{qu}}_{\mathbb{V},h,\beta,k,N} - M^{\mathrm{qu}}_{h,\beta,k,N}|$. By the definition of the median and an application of the Markov inequality, we have

$$\tfrac{1}{2} \leq \mathbb{P}[\log Z^{\mathrm{qu}}_{\mathbb{V},h,\beta,k,N} \geq M^{\mathrm{qu}}_{h,\beta,k,N}] \leq \mathbb{E} Z^{\mathrm{qu}}_{\mathbb{V},h,\beta,k,N} e^{-M^{\mathrm{qu}}_{h,\beta,k,N}}$$

and therefore

(1.9) $$M^{\mathrm{qu}}_{h,\beta,k,N} - \log \mathbb{E} Z^{\mathrm{qu}}_{\mathbb{V},h,\beta,k,N} \leq \log 2.$$

Now, let $t_{k,N} \stackrel{\mathrm{def}}{=} \log(\tfrac{1}{2} \mathbb{E} Z^{\mathrm{qu}}_{\mathbb{V},h,\beta,k,N}) - M^{\mathrm{qu}}_{h,\beta,k,N}$. Since we are looking for an upper bound, we can suppose $t_{k,N} \geq 0$. Again from (1.7), we obtain

(1.10)
$$\mathbb{P}\bigg[Z^{\mathrm{qu}}_{\mathbb{V},h,\beta,k,N} \geq \frac{1}{2} \mathbb{E} Z^{\mathrm{qu}}_{\mathbb{V},h,\beta,k,N}\bigg] = \mathbb{P}[\log Z^{\mathrm{qu}}_{\mathbb{V},h,\beta,k,N} - M^{\mathrm{qu}}_{h,\beta,k,N} \geq t_{k,N}]$$
$$\leq 4 \exp\bigg(\frac{-t_{k,N}^2}{16\beta^2 V^2 kN}\bigg).$$

Moreover, analogously to (1.6), we know that

(1.11) $$\mathbb{P}\bigg[Z^{\mathrm{qu}}_{\mathbb{V},h,\beta,k,N} \geq \frac{1}{2} \mathbb{E} Z^{\mathrm{qu}}_{\mathbb{V},h,\beta,k,N}\bigg] \geq \frac{1}{4} \frac{(\mathbb{E} Z^{\mathrm{qu}}_{\mathbb{V},h,\beta,k,N})^2}{\mathbb{E}(Z^{\mathrm{qu}}_{\mathbb{V},h,\beta,k,N})^2}.$$

For $\beta_0 < \beta_{\mathrm{c}}(h)$ chosen according to Theorem D, and $\varepsilon > 0$ and $k_\varepsilon$ according to Theorem E, we further have

(1.12) $$\frac{(\mathbb{E} Z^{\mathrm{qu}}_{\mathbb{V},h,\beta,k_\varepsilon,N})^2}{\mathbb{E}(Z^{\mathrm{qu}}_{\mathbb{V},h,\beta,k_\varepsilon,N})^2} \geq \varepsilon^2 \frac{(\mathbb{E} Z^{\mathrm{qu}}_{\mathbb{V},h,\beta,N})^2}{\mathbb{E}(Z^{\mathrm{qu}}_{\mathbb{V},h,\beta,N})^2} \geq \frac{\varepsilon^2}{K_{\mathrm{s.m.}}}$$

for all $\beta \leq \beta_0$ and $N \in \mathbb{N}$. A combination of (1.10), (1.11) and (1.12) then implies that

$$t_{k_\varepsilon,N} \leq \beta K_1 \sqrt{N}$$

for all $\beta \leq \beta_0$ and $N \in \mathbb{N}$, where the constant $K_1$ is given by

$$K_1 \stackrel{\mathrm{def}}{=} 4V \sqrt{k_\varepsilon \log \frac{16 K_{\mathrm{s.m.}}}{\varepsilon^2}}.$$

By the definition of $t_{k_\varepsilon,N}$, we therefore have

(1.13) $$\log \mathbb{E} Z^{\mathrm{qu}}_{\mathbb{V},h,\beta,k_\varepsilon,N} - M^{\mathrm{qu}}_{h,\beta,k_\varepsilon,N} \leq \beta K_1 \sqrt{N} + \log 2$$

and thus, as consequence of (1.8), (1.9) and (1.13),

(1.14) $$\mathbb{E}|\log \mathbb{E} Z^{\mathrm{qu}}_{\mathbb{V},h,\beta,k_\varepsilon,N} - \log Z^{\mathrm{qu}}_{\mathbb{V},h,\beta,k_\varepsilon,N}| \leq \beta K_2 \sqrt{N} + \log 2$$

for all $\beta \leq \beta_0$ and $N \in \mathbb{N}$, where $K_2 \stackrel{\mathrm{def}}{=} 8V\sqrt{k_\varepsilon \pi} + K_1$.

It remains to transfer (1.14) to the unmodified partition functions. By the triangle inequality, $|\log \mathbb{E} Z^{\mathrm{qu}}_{\mathbb{V},h,\beta,N} - \log Z^{\mathrm{qu}}_{\mathbb{V},h,\beta,N}|$ is bounded by

$$\log \mathbb{E} Z^{\mathrm{qu}}_{\mathbb{V},h,\beta,N} - \log \mathbb{E} Z^{\mathrm{qu}}_{\mathbb{V},h,\beta,k_\varepsilon,N} + |\log \mathbb{E} Z^{\mathrm{qu}}_{\mathbb{V},h,\beta,k_\varepsilon,N} - \log Z^{\mathrm{qu}}_{\mathbb{V},h,\beta,k_\varepsilon,N}|$$
$$- \log Z^{\mathrm{qu}}_{\mathbb{V},h,\beta,k_\varepsilon,N} + \log \mathbb{E} Z^{\mathrm{qu}}_{\mathbb{V},h,\beta,k_\varepsilon,N} - \log \mathbb{E} Z^{\mathrm{qu}}_{\mathbb{V},h,\beta,k_\varepsilon,N} + \log Z^{\mathrm{qu}}_{\mathbb{V},h,\beta,N}.$$

In order to handle the last summand in the above formula, recall also that

$$\mathbb{E} \log Z^{\mathrm{qu}}_{\mathbb{V},h,\beta,N} \leq \log \mathbb{E} Z^{\mathrm{qu}}_{\mathbb{V},h,\beta,N},$$

by Jensen's inequality. From (1.14) and Lemma E, we thus obtain

$$\mathbb{E}|\log \mathbb{E} Z^{\mathrm{qu}}_{\mathbb{V},h,\beta,N} - \log Z^{\mathrm{qu}}_{\mathbb{V},h,\beta,N}|$$
$$\leq 2 \log\left(\frac{\mathbb{E} Z^{\mathrm{qu}}_{\mathbb{V},h,\beta,N}}{\mathbb{E} Z^{\mathrm{qu}}_{\mathbb{V},h,\beta,k_\varepsilon,N}}\right) + 2\mathbb{E}|\log \mathbb{E} Z^{\mathrm{qu}}_{\mathbb{V},h,\beta,k_\varepsilon,N} - \log Z^{\mathrm{qu}}_{\mathbb{V},h,\beta,k_\varepsilon,N}|$$
$$\leq -2 \log \varepsilon + 2\beta K_2 \sqrt{N} + 2 \log 2$$

for all $\beta \leq \beta_0$ and $N \in \mathbb{N}$. By setting $K_{\mathrm{fr.e.}} \stackrel{\mathrm{def}}{=} 2\max\{K_2, \log 2 - \log \varepsilon\}$, the proof of Theorem C is completed. $\square$

For the proof of Theorem B, Theorem D and Lemma E, it remains to consider the annealed setting only. In Section 2, we deal with the "point-to-hyperplane" setting. That is, under the annealed path measure, we analyze finite paths $\mathcal{S}[n]$ with $S_1(n) = L$ for fixed $L \in \mathbb{N}$. In the first two subsections, such paths are approximated by more specific paths, the so-called bridges, and a renewal formalism is found by introducing irreducibility for bridges. In Section 2.3, we then prove the existence of a gap between the exponential decay rates of arbitrary and irreducible bridges.

In Section 3.1, we introduce an analogous renewal formalism for the "fixed-number-of-steps" setting. In Section 3.2, the exponential gap from Section 2.3 is transferred to that setting, implying the crucial part of Theorem B, namely analyticity of the Lyapunov exponent in the ballistic regime. With the renewal formalism and the exponential gap, we then also have the means to prove Lemma E, which is done in Section 3.3.

Finally, Section 4 is devoted to the proof of the second-moment estimate in Theorem D. It is based on a local decay estimate for sums of independent random variables and again on the gap between the exponential behaviors of arbitrary and irreducible bridges.

**2. Endpoint in given hyperplane.** As in the previous section, we consider a nearest-neighbor random walk $\mathcal{S} = (S(n))_{n \in \mathbb{N}_0}$ on $\mathbb{Z}^d$, with start at the origin and drift $h$ in the direction of the first axis, defined on a probability space $(\Omega, \mathcal{F}, P_h)$.

In this section, we investigate the behavior of finite random paths with start at the origin and endpoint in a hyperplane

$$\mathcal{H}_L \stackrel{\text{def}}{=} \{(\xi_1, \ldots, \xi_d) \in \mathbb{Z}^d : \xi_1 = L\},$$

exponentially weighted by a nonrandom path potential $\Phi_\beta$. We therefore allow a more general setting for $\Phi_\beta$ than for the annealed path potential $\Phi_\beta^{\text{an}}$ from Section 1. More precisely, we assume that $\varphi : \mathbb{R}^+ \to \mathbb{R}^+$ is a nonconstant, concave increasing function with

(2.1) $$\lim_{t \to 0} \varphi(t) = \varphi(0) = 0$$

and

(2.2) $$\lim_{t \to \infty} \varphi(t) = \infty,$$

(2.3) $$\lim_{t \to \infty} \frac{\varphi(t)}{t} = 0.$$

For $\beta \geq 0$, we define $\varphi_\beta : \mathbb{R}^+ \to \mathbb{R}^+$ by

$$\varphi_\beta(t) \stackrel{\text{def}}{=} \varphi(t\beta), \qquad t \in \mathbb{R}^+.$$

The function $\varphi_\beta$ plays the role of the annealed potential $\varphi_\beta^{\text{an}}$ from Section 1, with coincidence in the case

(2.4) $$\varphi(t) = -\log \mathbb{E} \exp(-tV_x).$$

Assumption (2.3) is needed for Theorem A only. Under (2.4), it corresponds to the assumption $\operatorname{ess\,inf} V_x = 0$ from the quenched setting, and (2.2) is equivalent to $P_h[V_x = 0] < 1$.

The path potential is now defined as in the annealed setting from the previous section: for $\beta \geq 0$ and $N, M \in \mathbb{N}_0$ with $M \leq N$, set

$$\Phi_\beta(N) \stackrel{\text{def}}{=} \sum_{x \in \mathbb{Z}^d} \varphi_\beta(\ell_x(N)),$$

$$\Phi_\beta(M, N) \stackrel{\text{def}}{=} \sum_{x \in \mathbb{Z}^d} \varphi_\beta(\ell_x(M, N)),$$

where

$$\ell_x(N) \stackrel{\text{def}}{=} \sum_{n=1}^N 1_{\{S(n)=x\}}, \qquad \ell_x(M, N) \stackrel{\text{def}}{=} \sum_{n=M+1}^N 1_{\{S(n)=x\}}$$

denote the number of visits to the site $x \in \mathbb{Z}^d$ by the random walk $\mathcal{S}[1, N]$, respectively $\mathcal{S}[M+1, N]$.

We derive some elementary properties of the path potential $\Phi_\beta$ arising from the assumptions on $\varphi_\beta$. For $N, M \in \mathbb{N}_0$ with $M \leq N$, let

$$R(N) \stackrel{\text{def}}{=} \{x \in \mathbb{Z}^d : \ell_x(N) \geq 1\},$$

$$R(M, N) \stackrel{\text{def}}{=} \{x \in \mathbb{Z}^d : \ell_x(M, N) \geq 1\}$$

denote the sets of sites visited by $\mathcal{S}[1, N]$, respectively $\mathcal{S}[M+1, N]$.

LEMMA 2.1. (a) *For any $\beta \geq 0$ and $N \in \mathbb{N}_0$, we have*

$$\varphi_\beta(1) \sharp R(N) \leq \Phi_\beta(N) \leq \varphi_\beta(1) N.$$

(b) *For any $\beta \geq 0$ and $M, N \in \mathbb{N}_0$ with $M \leq N$, we have*

$$\Phi_\beta(N) \leq \Phi_\beta(M) + \Phi_\beta(M, N).$$

*Moreover, if $\omega \in \{R(M) \cap R(M, N) = \varnothing\}$, we have*

$$\Phi_\beta(N, \omega) = \Phi_\beta(M, \omega) + \Phi_\beta(M, N, \omega).$$

(c) *For any $\beta \geq 0$ and $M_1, M_2, N \in \mathbb{N}_0$ with $M_1 \leq M_2 \leq N$, we have*

$$\Phi_\beta(N) \geq \Phi_\beta(M_1, M_2).$$

*Moreover, if $\omega \in \{R(M_1) \cap R(M_2, N) = \varnothing\}$, we have*

$$\Phi_\beta(N, \omega) \geq \Phi_\beta(M_1, \omega) + \Phi_\beta(M_2, N, \omega).$$

PROOF. (a) For the lower bound, we use the monotonicity of $\varphi_\beta$ to obtain

$$\Phi_\beta(N) = \sum_{x \in R(N)} \varphi_\beta(\ell_x(N)) \geq \sum_{x \in R(N)} \varphi_\beta(1) = \varphi_\beta(1) \sharp R(N).$$

For the upper estimate, we inductively apply the concavity of $\varphi_\beta$ to obtain

$$\Phi_\beta(N) = \sum_{x \in \mathbb{Z}^d} \varphi_\beta(\ell_x(N)) \leq \sum_{x \in \mathbb{Z}^d} \varphi_\beta(1) \ell_x(N) = \varphi_\beta(1) N.$$

(b) By the concavity of $\varphi_\beta$ and (2.1), we have

$$\Phi_\beta(N) = \sum_{x \in \mathbb{Z}^d} \varphi_\beta(\ell_x(M) + \ell_x(M, N))$$

$$\leq \sum_{x \in \mathbb{Z}^d} \varphi_\beta(\ell_x(M)) + \varphi_\beta(\ell_x(M, N))$$

$$= \Phi_\beta(M) + \Phi_\beta(M, N).$$

For $\omega \in \{R(M) \cap R(M,N) = \varnothing\}$, we have

$$\Phi_\beta(N,\omega) = \sum_{x \in R(M,\omega)} \varphi_\beta(\ell_x(N,\omega)) + \sum_{x \in R(M,N,\omega)} \varphi_\beta(\ell_x(N,\omega))$$

$$= \sum_{x \in R(M,\omega)} \varphi_\beta(\ell_x(M,\omega)) + \sum_{x \in R(M,N,\omega)} \varphi_\beta(\ell_x(M,N,\omega))$$

$$= \Phi_\beta(M,\omega) + \Phi_\beta(M,N,\omega).$$

(c) The monotony of $\varphi_\beta$ implies that

$$\Phi_\beta(N) = \sum_{x \in \mathbb{Z}^d} \varphi_\beta(\ell_x(N)) \geq \sum_{x \in \mathbb{Z}^d} \varphi_\beta(\ell_x(M_1, M_2)) = \Phi_\beta(M_1, M_2).$$

For $\omega \in \{R(M_1) \cap R(M_2, N) = \varnothing\}$, we have

$$\Phi_\beta(N,\omega) \geq \sum_{x \in R(M_1,\omega)} \varphi_\beta(\ell_x(N,\omega)) + \sum_{x \in R(M_2,N,\omega)} \varphi_\beta(\ell_x(N,\omega))$$

$$\geq \sum_{x \in R(M_1,\omega)} \varphi_\beta(\ell_x(M_1,\omega)) + \sum_{x \in R(M_2,N,\omega)} \varphi_\beta(\ell_x(M_2,N,\omega))$$

$$= \Phi_\beta(M_1,\omega) + \Phi_\beta(M_2,N,\omega). \qquad \square$$

2.1. *Masses for paths and bridges.* We start with a few comments on the process of the first components of $\mathcal{S}$, that is,

$$\mathcal{S}_1 \stackrel{\text{def}}{=} (S_1(n))_{n \in \mathbb{N}_0}.$$

The process $\mathcal{S}_1$ is itself a random walk on $\mathbb{Z}$, again with independent increments and drift in the positive direction. It can be expressed by

$$S_1(n) = \sum_{m=1}^{n} (S_1(m) - S_1(m-1))$$

for $n \in \mathbb{N}$, where the random variables $(S_1(m) - S_1(m-1))_{m \in \mathbb{N}}$ are independent and identically distributed. Since $E_h S_1(n) > 0$ for $h > 0$, we then have

$$P_h[S_1(n) \to \infty \text{ as } n \to \infty] = 1,$$

by the strong law of large numbers. This convergence property, as we show next, implies transience to the process $\mathcal{S}_1$, which is here equivalent to the fact that the probability

(2.5) $$\alpha(h) \stackrel{\text{def}}{=} P_h[S_1(n) > 0 \text{ for all } n \in \mathbb{N}]$$

is strictly greater than zero: for $h > 0$, and with $m \in \mathbb{N}_0$ denoting the last time the random walk $\mathcal{S}_1$ is in $L \in \mathbb{N}_0$, we have

$$1 = \sum_{m=0}^{\infty} P_h[S_1(m) = L, S_1(m+n) > L \text{ for all } n \in \mathbb{N}]$$

$$= \sum_{m=0}^{\infty} P_h[S_1(n) = L] P_h[S_1(n) > 0 \text{ for all } n \in \mathbb{N}]$$

and therefore

$$(2.6) \qquad \sum_{m=0}^{\infty} P_h[S_1(m) = L] = \frac{1}{\alpha(h)} < \infty.$$

Observe, in particular, that the left-hand side of (2.6) does not depend on $L \in \mathbb{N}_0$.

REMARK 2.2. For every $h > 0$ and $L \in \mathbb{N}$, we have

$$(2.7) \qquad P_h[H_{-L} < \infty] = e^{-2Lh},$$

where the stopping time

$$H_{-L} \stackrel{\text{def}}{=} \inf\{n \in \mathbb{N} : S_1(n) = -L\}$$

denotes the time of the random walk's first visit to the hyperplane $\mathcal{H}_{-L}$.

PROOF. The Markov property implies that

$$P_h[H_{-L} < \infty] = P_h[H_{-1} < \infty]^L.$$

We can thus restrict our attention to the case $L = 1$. Also by the Markov property, we have

$$P_h[H_{-1} < \infty] = P_h[S_1(1) = -1] + P_h[S_1(1) = 0] P_h[H_{-1} < \infty]$$
$$+ P_h[S_1(1) = 1] P_h[H_{-1} < \infty]^2,$$

which is a quadratic equation in the variable $P_h[H_{-1} < \infty]$, with solutions $1$ and $P_h[S_1(1) = -1]/P[S_1(1) = 1]$. To find the correct one among these two solutions, observe that

$$\alpha(h) = P_h[S_1(1) = 1] P_h[H_{-1} = \infty],$$

again by the Markov property. We thus have

$$P_h[H_{-1} < \infty] = 1 - \frac{\alpha(h)}{P_h[S_1(1) = 1]} < 1$$

since $\alpha(h) > 0$ by (2.6), and therefore

$$P_h[H_{-1} < \infty] = \frac{P_h[S_1(1) = -1]}{P[S_1(1) = 1]} = e^{-2h},$$

where, at the second step, the concrete definition of $P_h$ in Section 1 is used. □

We now return to random walks on $\mathbb{Z}^d$. In the present setting, for any $h > 0$, $\beta \geq 0$ and $L \in \mathbb{N}$, the counterpart of $\overline{Z}^{\text{an}}_{h,\beta,L}$ from the point-to-hyperplane setting of Section 1 is given by

(2.8) $$\overline{G}_{h,\beta}(L) \stackrel{\text{def}}{=} \sum_{N=1}^{\infty} E_h[\exp(-\Phi_\beta(N)); \{S_1(N) = L\}].$$

REMARK 2.3. The drift $h > 0$ ensures that (2.8) is finite. More exactly, for any $h > 0$, $\beta \geq 0$ and $L \in \mathbb{N}$, we have

$$P_h[S(1) = 1]^L e^{-\varphi_\beta(1)L} \leq \overline{G}_{h,\beta}(L) \leq \frac{1}{\alpha(h)} e^{-\varphi_\beta(1)L}.$$

PROOF. It is plain that

$$\overline{G}_{h,\beta}(L) \geq E_h[\exp(-\Phi_\beta(L)); \{S_1(L) = L\}],$$

which implies the lower estimate. For $\omega \in \{S_1(N) = L\}$, we obviously have $L \leq \sharp R(N, \omega)$. From Lemma 2.1(a), we thus obtain

$$E_h[\exp(-\Phi_\beta(L)); \{S_1(N) = L\}] \leq e^{-\varphi_\beta(1)L} P_h[S_1(N) = L]$$

for all $L \in \mathbb{N}$, from which the lower estimate now follows by (2.6). □

We are interested in the exponential behavior of $\overline{G}_{h,\beta}(L)$ as $L \to \infty$. This behavior is easier to study when the expectations in (2.8) are restricted to so-called bridges.

DEFINITION 2.4. Suppose $\omega \in \Omega$ and $N, M \in \mathbb{N}_0$ with $M \leq N$. The finite path $\mathcal{S}[M,N](\omega)$ is called a *bridge* if

$$S_1(M, \omega) < S_1(n, \omega) \leq S_1(N, \omega)$$

is valid for $n = M+1, \ldots, N$. In that case, the *span* of the bridge $\mathcal{S}[M,N](\omega)$ is given by $S_1(N, \omega) - S_1(M, \omega)$.

For $L \in \mathbb{N}$ and $M, N \in \mathbb{N}_0$ with $M \leq N$, we define

$$\mathrm{br}(L; N) \stackrel{\mathrm{def}}{=} \{\mathcal{S}[N] \text{ is a bridge of span } L\},$$

$$\mathrm{br}(L; M, N) \stackrel{\mathrm{def}}{=} \{\mathcal{S}[M, N] \text{ is a bridge of span } L\}.$$

For $h, \beta \geq 0$ and $L \in \mathbb{N}$, we further define

$$b_{h,\beta}(L; N) \stackrel{\mathrm{def}}{=} E_h[\exp(-\Phi_\beta(N)); \mathrm{br}(L; N)],$$

$$\overline{B}_{h,\beta}(L) \stackrel{\mathrm{def}}{=} \sum_{N=1}^{\infty} b_{h,\beta}(L).$$

REMARK 2.5. For any $h, \beta \geq 0$ and $L \in \mathbb{N}$, we have

$$P_h[S(1) = 1]^L e^{-\varphi_\beta(1)L} \leq \overline{B}_{h,\beta}(L) \leq e^{-\varphi_\beta(1)L}.$$

PROOF. The lower estimate is proved as in Remark 2.3. For the upper estimate, observe that

$$b_{h,\beta}(L; N) \leq b_{h,0}(L; N) e^{-\varphi_\beta(1)L}$$

for $L \leq N$ by Lemma 2.1(a). By the Markov property, we furthermore have

$$b_{h,0}(L; N) = P_h[0 \leq S_1(n) < S_1(N) = L \text{ for } 0 \leq n < N]$$

and thus

$$\overline{B}_{h,\beta}(L) \leq \sum_{N=1}^{\infty} P_h[H_L = N] e^{-\varphi_\beta(1)L} = P_h[H_L < \infty] e^{-\varphi_\beta(1)L},$$

where $H_L \stackrel{\mathrm{def}}{=} \min\{n \in \mathbb{N} : S_1(n) = L\}$. □

The well-known subadditive limit lemma (see, e.g., page 9 of [11]) states the following. Let $(a_n)_{n \in \mathbb{N}}$ be a sequence of real numbers with the *subadditivity property*

$$a_{n+m} \leq a_n + a_m$$

for all $m, n \in \mathbb{N}$. We then have

$$\lim_{n \to \infty} \frac{a_n}{n} = \inf\left\{\frac{a_n}{n} : n \in \mathbb{N}\right\}.$$

We want to apply the subadditive limit lemma to $(-\log \overline{B}_{h,\beta})_{L \in \mathbb{N}}$. The subadditivity property is a consequence of the following lemma.

LEMMA 2.6. *For $h, \beta \geq 0$ and $L_1, L_2 \in \mathbb{N}$, we have*

$$\overline{B}_{h,\beta}(L_1 + L_2) \geq \overline{B}_{h,\beta}(L_1)\overline{B}_{h,\beta}(L_2).$$

*Moreover, for nonvanishing drift $h > 0$, we have*

$$\overline{B}_{h,\beta}(L_1 + L_2) \leq \frac{1}{\alpha(h)}\overline{B}_{h,\beta}(L_1)\overline{B}_{h,\beta}(L_2).$$

PROOF. In order to obtain the lower estimate, observe that

$$\mathrm{br}(L_1 + L_2; N) \supset \bigcup_{M=1}^{N-1} \mathrm{br}(L_1; M) \cap \mathrm{br}(L_2; M, N),$$

where the right-hand side is a union of disjoint sets and where $M$ denotes the time of the unique visit of $\mathcal{S}[N]$ to the hyperplane $\mathcal{H}_{L_1}$. For $\omega \in \mathrm{br}(L_1; M) \cap \mathrm{br}(L_2; M, N)$, we further have

$$\Phi_\beta(N, \omega) = \Phi_\beta(M, \omega) + \Phi_\beta(M, N, \omega)$$

by Lemma 2.1(b). By splitting over all possible values of $M$ and using the Markov property to renew the random walk $\mathcal{S}$ at that time, we obtain

$$\overline{B}_{h,\beta}(L_1 + L_2) \geq \sum_{N=1}^{\infty} \sum_{M=1}^{N-1} E_h[e^{-\Phi_\beta(M)}1_{\mathrm{br}(L_1;M)}e^{-\Phi_\beta(M,N)}1_{\mathrm{br}(L_2;M,N)}]$$

$$= \sum_{N=1}^{\infty} \sum_{M=1}^{N-1} b_{h,\beta}(L_1; M)b_{h,\beta}(L_2; N - M)$$

$$= \overline{B}_{h,\beta}(L_1)\overline{B}_{h,\beta}(L_2),$$

which proves the lower estimate of the lemma.

The upper estimate is shown in a similar way. The event $\mathrm{br}(L_1 + L_2; N)$ is contained in the union

$$\bigcup_{M_1=1}^{N-1} \bigcup_{M_2=M_1}^{N-1} \mathrm{br}(L_1; M_1) \cap \{S_1(M_2) = L_1\} \cap \mathrm{br}(L_2; M_2, N),$$

in which $M_1$, for $\omega \in \mathrm{br}(L_1 + L_2; N)$, may be chosen as the time of the first visit of $\mathcal{S}[N](\omega)$ to the hyperplane $\mathcal{H}_{L_1}$, and where $M_2$ stands for the time of the last visit to $\mathcal{H}_{L_1}$. By splitting over all possible values of $M_1$ and $M_2$, and using Lemma 2.1(c) and the Markov property to renew the random walk at these times, we obtain

$$\overline{B}_{h,\beta}(L_1 + L_2)$$

$$\leq \sum_{N=1}^{\infty} \sum_{M_1=1}^{N-1} \sum_{M_2=M_1}^{N-1} E_h[e^{-\Phi_\beta(M_1)}1_{\mathrm{br}(L_1;M_1)}1_{\{S_1(M_2)=L_1\}}$$

$$\times e^{\Phi_\beta(M_2,N)} 1_{\mathrm{br}(L_2;M_2,N)}]$$

$$= \sum_{N=1}^{\infty} \sum_{M_1=1}^{N-1} \sum_{M_2=M_1}^{N-1} b_{h,\beta}(L_1;M_1) P_h[S_1(M_2-M_1)=0] b_{h,\beta}(L_2;N-M_2)$$

$$= \overline{B}_{h,\beta}(L_1) \overline{B}_{h,\beta}(L_2) \sum_{m=0}^{\infty} P_h[S_1(m)=0],$$

from which the upper estimate of the lemma follows by (2.6). □

PROPOSITION 2.7. *For any $h, \beta \geq 0$, the mass*

$$\overline{m}_B(h,\beta) \stackrel{\mathrm{def}}{=} \lim_{L \to \infty} \frac{-\log \overline{B}_{h,\beta}(L)}{L}$$

*of $\overline{B}_{h,\beta}$ exists in $[\varphi_\beta(1), \infty)$, is continuous as function on $\mathbb{R}^+ \times \mathbb{R}^+$ and satisfies*

(2.9) $$\overline{B}_{h,\beta}(L) \leq e^{-\overline{m}_B(h,\beta)L}$$

*for all $L \in \mathbb{N}$. Moreover, for nonvanishing $h > 0$, we have*

(2.10) $$\overline{B}_{h,\beta}(L) \geq \alpha(h) e^{-\overline{m}_B(h,\beta)L}$$

*for all $L \in \mathbb{N}$.*

PROOF. The sequence $(-\log \overline{B}_{h,\beta}(L))_{L \in \mathbb{N}}$ is subadditive by Lemma 2.6. The subadditive limit lemma thus yields

$$\overline{m}_B(h,\beta) = \inf\left\{\frac{-\log \overline{B}_{h,\beta}(L)}{L} : L \in \mathbb{N}\right\} \in [-\infty, \infty),$$

which includes the existence of the limit and implies the estimate in (2.9). The lower bound $\varphi_\beta(1)$ for the mass $\overline{m}_B(h,\beta)$ follows from Remark 2.3.

By the above expression, as an infimum of continuous functions, the mass $\overline{m}_B$ is upper semicontinuous. In order to obtain lower semicontinuity, it is convenient to consider

$$\widehat{B}_{\lambda,\beta}(L) \stackrel{\mathrm{def}}{=} \sum_{N=1}^{\infty} E_0[\exp(-\Phi_\beta(N) - \lambda N); \mathrm{br}(L;N)]$$

for $\lambda \geq 0$ and $L \in \mathbb{N}$. By the definition of $P_h$ in (1.1), we have

$$\widehat{m}_B(\lambda_h,\beta) \stackrel{\mathrm{def}}{=} \lim_{L \to \infty} \frac{-\log \widehat{B}_{\lambda_h,\beta}(L)}{L} = \overline{m}_B(h,\beta) + h,$$

where $\lambda_h = \log E_0[\exp(h \cdot S_1(1))]$. It consequently suffices to show that $\widehat{m}_B$ is lower semicontinuous.

To see this, observe that for any fixed $N \in \mathbb{N}$, the map $\beta \mapsto \Phi_\beta(N)$ inherits the concavity of $\varphi$. By the Hölder inequality, for any $(\lambda, \beta), (\lambda', \beta') \in \mathbb{R}^+ \times \mathbb{R}^+$ and $t \in [0, 1]$, we thus have

$$\widehat{B}_{t\lambda+(1-t)\lambda', t\beta+(1-t)\beta'}(L)$$
$$\leq \sum_{N=1}^{\infty} E_0[e^{-t\Phi_\beta(N)-(1-t)\Phi_{\beta'}(N)-t\lambda-(1-t)\lambda'}; \mathrm{br}(L; N)]$$
$$\leq \sum_{N=1}^{\infty} E_0[e^{-\Phi_\beta(N)-\lambda N}; \mathrm{br}(L; N)]^t E_0[e^{-\Phi_{\beta'}(N)-\lambda' N}; \mathrm{br}(L; N)]^{(1-t)}$$
$$\leq \widehat{B}_{\lambda,\beta}(L)^t \widehat{B}_{\lambda',\beta'}(L)^{(1-t)}.$$

Therefore, for any fixed $L \in \mathbb{N}$, the negative logarithm of $\widehat{B}_{\lambda,\beta}(L)$ is concave as a function of $(\lambda, \beta)$. The mass $\widehat{m}_B$ inherits this concavity and, as a consequence, is lower semicontinuous.

It remains to prove (2.10). To this end, we consider

$$\widetilde{B}_{h,\beta}(L) \stackrel{\mathrm{def}}{=} \alpha(h)^{-1} \overline{B}_{h,\beta}(L), \qquad L \in \mathbb{N}.$$

The sequence $(\log \widetilde{B}_{h,\beta}(L))_{L \in \mathbb{N}}$ is subadditive by Lemma 2.6. As a consequence, we have

$$\alpha(h)^{-1} \overline{B}_{h,\beta}(L) = \widetilde{B}_{h,\beta}(L) \geq e^{-\widetilde{m}_B(h,\beta)L}$$

for all $L \in \mathbb{N}$, by the subadditivity limit lemma. Thereby, $\widetilde{m}_B(h, \beta)$ denotes the mass of $\widetilde{B}_{h,\beta}$ and is given by

$$\widetilde{m}_B(h, \beta) \stackrel{\mathrm{def}}{=} \lim_{L \to \infty} \frac{-\log \widetilde{B}_{h,\beta}(L)}{L}$$
$$= \lim_{L \to \infty} \frac{-\log \overline{B}_{h,\beta}(L)}{L} + \frac{\log \alpha(h)}{L} = \overline{m}_B(h, \beta). \qquad \square$$

By means of the following lemma, the results on $\overline{B}_{h,\beta}$ in Proposition 2.7 can be transfered to $\overline{G}_{h,\beta}$.

LEMMA 2.8. *For any $h > 0$ and $\beta \geq 0$, we have*
$$\alpha(h)^2 \overline{G}_{h,\beta}(L) \leq \overline{B}_{h,\beta}(L)$$
*for all $L \in \mathbb{N}$.*

PROOF. The event $\{S_1(N) = L\}$ is contained in the union

$$\bigcup_{M_1=0}^{N-1} \bigcup_{M_2=M_1+1}^{N} \{S_1(M_1) = 0\} \cap \mathrm{br}(L; M_1, M_2) \cap \{S_1(N) = L\},$$

in which $M_1$ and $M_2$, for $\omega \in \{S_1(N) = L\}$, may be chosen as the time of the last return of $\mathcal{S}[N](\omega)$ to the hyperplane $\mathcal{H}_0$, respectively the first visit of $\mathcal{S}[M_1, N](\omega)$ to the hyperplane $\mathcal{H}_L$. By splitting over all possible values of $M_1$ and $M_2$, and applying Lemma 2.1(c) and the Markov property to renew the random walk at these times, we obtain that $\overline{G}_{h,\beta}(L)$ is bounded by

$$\sum_{N=1}^{\infty} \sum_{M_1=0}^{N-1} \sum_{M_2=M_1+1}^{N} E_h[1_{\{S_1(M_1)=0\}} e^{-\Phi_\beta(M_1, M_2)} 1_{\mathrm{br}(L; M_1, M_2)} 1_{\{S_1(N)=L\}}]$$

$$= \sum_{N=1}^{\infty} \sum_{M_1=0}^{N-1} \sum_{M_2=M_1+1}^{N} P_h[S_1(M_1) = 0] b_{h,\beta}(L; M_2 - M_1)$$

$$\times P_h[S_1(N - M_2) = 0]$$

$$= \overline{B}_{h,\beta}(L) \left( \sum_{n=0}^{\infty} P_h[S_1(n) = 0] \right)^2,$$

from which the lemma follows, by (2.6). $\square$

COROLLARY 2.9. *For any $h > 0$ and $\beta \geq 0$, we have*

$$\overline{m}_G(h, \beta) \stackrel{\text{def}}{=} \lim_{L \to \infty} \frac{-\log \overline{G}_{h,\beta}(L)}{L} = \overline{m}_B(h, \beta)$$

*and*

(2.11) $$\alpha(h) e^{-\overline{m}_G(h,\beta)L} \leq \overline{G}_{h,\beta}(L) \leq \frac{1}{\alpha(h)^2} e^{-\overline{m}_G(h,\beta)L}$$

*for all $L \in \mathbb{N}$.*

PROOF. The corollary follows from Proposition 2.7 and Lemma 2.8. $\square$

REMARK 2.10. For any $h > 0$, as anticipated in Section 1.3, there exists a unique parameter $\beta_c(h) > 0$, such that

$$\overline{m}_G(0, \beta_c(h)) = h.$$

PROOF. For any $h > 0$, by Corollary 2.9 and as shown in the proof of Proposition 2.7, the mass $\overline{m}_G(h, \beta)$ is continuous and concave increasing in the variable $\beta \in \mathbb{R}^+$. Furthermore, we have $\overline{m}_G(h, 0) = 0$ by (2.6) and the assumption $\varphi(0) = 0$, and $\lim_{\beta \to \infty} \overline{m}_G(h, \beta) = \infty$ by Remark 2.3 and the assumption $\lim_{t \to \infty} \varphi(t) = \infty$. This limiting behavior for $\beta \to \infty$, in combination with the concavity of the mass shown above, moreover yields that the monotonicity of $\overline{m}_G(h, \beta)$ in $\beta \in \mathbb{R}^+$ is strict. This completes the proof of the remark. $\square$

2.2. *Irreducible bridges and renewal results.* For the rest of the this section, we fix an integer $p \in \mathbb{N}$ and consider independent copies

$$\mathcal{S}^j = (S^j(n))_{n \in \mathbb{N}_0}, \qquad j = 1, \ldots, p,$$

of the random walk $\mathcal{S}$. We assume these copies to be defined on the product space $(\Omega^p, \mathcal{F}^{\otimes p})$, on which the $p$-fold product measure, in order to keep notation simple, is again denoted by $P_h$. We compose from $\mathcal{S}^1, \ldots, \mathcal{S}^p$ a random process $\mathcal{S}^{(p)}$ with values in $(\mathbb{Z}^d)^p$ by setting

$$S^{(p)}(n) \stackrel{\text{def}}{=} (S^1(n^1), \ldots, S^p(n^p)),$$

for $n = (n^1, \ldots, n^p) \in \mathbb{N}_0^p$, and

$$\mathcal{S}^{(p)} \stackrel{\text{def}}{=} (S^{(p)}(n))_{n \in \mathbb{N}_0^p}.$$

The process $\mathcal{S}^{(p)}$ inherits the Markov property from $\mathcal{S}^1, \ldots, \mathcal{S}^p$ in the following way: for $M = (M^1, \ldots, M^p)$ and $N = (N^1, \ldots, N^p) \in \mathbb{N}_0^p$, we write

$$M \leq N, \quad \text{if and only if} \quad M^j \leq N^j \quad \text{for } j = 1, \ldots, p,$$
$$M < N, \quad \text{if and only if} \quad M^j < N^j \quad \text{for } j = 1, \ldots, p$$

and set

$$[M, \ldots, N] \stackrel{\text{def}}{=} \{n \in \mathbb{N}^p : M \leq n \leq N\}.$$

The origin in $\mathbb{N}_0^p$ is denoted by $0$. Suppose now that $M, N \in \mathbb{N}_0^p$ with $M \leq N$, and $x_n = (x_{n^1}^1, \ldots, x_{n^p}^p) \in (\mathbb{Z}^d)^p$ for $n = (n^1, \ldots, n^p) \in [0, \ldots, N]$. Then, if

$$P_h[S^{(p)}(m) = x_m \text{ for } m \in [0, \ldots, M]] > 0$$

is valid, we have

$$P_h[S^{(p)}(n) = x_n \text{ for } n \in [M, \ldots, N] | S^{(p)}(m) = x_m \text{ for } m \in [0, \ldots, M]]$$
$$= \prod_{j=1}^p P_h[\mathcal{S}^j[M^j, N^j] = (x_{M^j}^j, \ldots, x_{N^j}^j) | \mathcal{S}^j[M^j] = (x_0^j, \ldots, x_{M^j}^j)]$$
$$= \prod_{j=1}^p P_h[\mathcal{S}^j[N^j - M^j] = (x_{M^j}^j - x_{M^j}^j, \ldots, x_{N^j}^j - x_{M^j}^j)]$$
$$= P_h[S^{(p)}(n) = x_{M+n} - x_M \text{ for } n \in [0, \ldots, N - M]].$$

This means that, similarly to $\mathcal{S}$ in the previous subsection, the process $\mathcal{S}^{(p)}$ can be renewed at any time $M$.

In Definition 2.4, the denomination bridge was introduced in the context of a single random walk. We want to generalize it to the present setting: for $M, N \in \mathbb{N}_0^p$ with $M \leq N$, we write

$$\mathcal{S}^{(p)}[N] \stackrel{\text{def}}{=} (S^{(p)}(n))_{n \in \mathbb{N}_0^p : n \leq N},$$

$$\mathcal{S}^{(p)}[M, N] \stackrel{\text{def}}{=} (S^{(p)}(n))_{n \in \mathbb{N}_0^p : M \leq n \leq N}$$

for finite subpaths of $\mathcal{S}^{(p)}$ in $(\mathbb{Z}^d)^p$.

DEFINITION 2.11. Suppose $\omega \in \Omega^p$ and $M, N \in \mathbb{N}_0^p$ with $M \leq N$. The finite path $\mathcal{S}^{(p)}[M, N](\omega)$ is called a *bridge* if

$$S_1^1(M^1, \omega) = S_1^j(M^j, \omega) < S_1^j(n^j, \omega) \leq S_1^j(N^j, \omega) = S_1^1(N^1, \omega)$$

is valid for all $n^j = M^j + 1, \ldots, N^j$ and $j = 1, \ldots, p$. In that case, the *span* of the bridge $\mathcal{S}^{(p)}[M, N](\omega)$ is given by $S_1^1(N^1, \omega) - S_1^1(M^1, \omega)$.

REMARK. By this definition, a finite path $\mathcal{S}^{(p)}[M, N]$ in $(\mathbb{Z}^d)^p$ is a bridge if and only if $\mathcal{S}^1[M^1, N^1], \ldots, \mathcal{S}^p[M^p, N^p]$ are bridges in the sense of Definition 2.4, with start and endpoint each in a common hyperplane. Also, observe that the definition includes the case $M = N$, in which $\mathcal{S}^{(p)}[M, N]$ is a bridge of span zero.

At the beginning of this section, we introduced the path potential $\Phi_\beta$ for the random walk $\mathcal{S}$. For any $j \in \{1, \ldots, p\}$, we denote the corresponding potential associated with $\mathcal{S}^j$ by $\Phi_\beta^j$. For $M, N \in \mathbb{N}_0^p$ with $M \leq N$, a path potential for the process $\mathcal{S}^{(p)}$ is then given by

$$\Phi_\beta^{(p)}(N) \stackrel{\text{def}}{=} \sum_{j=1}^p \Phi_\beta^j(N^j),$$

(2.12)

$$\Phi_\beta^{(p)}(M, N) \stackrel{\text{def}}{=} \sum_{j=1}^p \Phi_\beta^j(M^j, N^j).$$

For $h, \beta \geq 0$, $L \in \mathbb{N}_0$ and $M, N \in \mathbb{N}_0^p$ with $M \leq N$, we now define

$$\mathrm{br}^p(L; N) \stackrel{\text{def}}{=} \{\mathcal{S}^{(p)}[N] \text{ is a bridge of span } L\},$$

$$\mathrm{br}^p(L; M, N) \stackrel{\text{def}}{=} \{\mathcal{S}^{(p)}[M, N] \text{ is a bridge of span } L\}$$

and also let

$$b_{h,\beta}^p(L; N) \stackrel{\text{def}}{=} E_h[\exp(-\Phi_\beta^{(p)}(N)); \mathrm{br}^p(L; N)],$$

$$\overline{B}_{h,\beta}^p(L) \stackrel{\text{def}}{=} \sum_{N \in \mathbb{N}_0^p} b_{h,\beta}^p(L; N).$$

REMARK 2.12. By the independence of $\mathcal{S}^1, \ldots, \mathcal{S}^p$, we have

$$\overline{B}^p_{h,\beta}(L) = (\overline{B}_{h,\beta}(L))^p \tag{2.13}$$

for all $L \in \mathbb{N}$. As a consequence, the mass

$$\overline{m}^p_B(h, \beta) \stackrel{\text{def}}{=} \lim_{L \to \infty} \frac{-\log \overline{B}^p_{h,\beta}(L)}{L}$$

of $\overline{B}^p_{h,\beta}$ exists and $\overline{m}^p_B(h, \beta) = p\overline{m}_B(h, \beta)$. By Proposition 2.7, we further have

$$\overline{B}^p_{h,\beta}(L) \leq e^{-\overline{m}^p_B(h,\beta)L} \tag{2.14}$$

for all $L \in \mathbb{N}$. Moreover, for nonvanishing $h > 0$, we have

$$\overline{B}^p_{h,\beta}(L) \geq \alpha(h)^p e^{-\overline{m}^p_B(h,\beta)L} \tag{2.15}$$

for all $L \in \mathbb{N}$.

Bridges allow a treatment using the tools of renewal theory. The decisive concepts for this treatment are the following.

DEFINITION 2.13. Suppose that $L \in \mathbb{N}$ and $\omega \in \text{br}^p(L; M, N)$ for $M \in \mathbb{N}_0^p$ and $N \in \mathbb{N}^p$ with $M < N$. An integer $R$ with $S_1^1(M^1, \omega) < R \leq S_1^1(N^1, \omega)$ is called a *breaking point* of $\mathcal{S}^{(p)}[M, N](\omega)$ if there exists $n \in \mathbb{N}^p$ with $M < n \leq N$ such that

   (i) $\mathcal{S}^{(p)}[M, n](\omega)$ is a bridge of span $R - S_1^1(M^1)$,
   (ii) $\mathcal{S}^{(p)}[n, N](\omega)$ is a bridge of span $S_1^1(N^1) - R$.

Moreover, the bridge $\mathcal{S}^{(p)}[M, N](\omega)$ is called *irreducible* if $S_1^1(N^1, \omega)$ is its only breaking point.

For $h, \beta \geq 0$, $L \in \mathbb{N}$, $M \in \mathbb{N}_0^p$ and $N \in \mathbb{N}^p$ with $M < N$, we now set

$$\text{ir}^p(L; N) \stackrel{\text{def}}{=} \{\mathcal{S}^{(p)}[N] \text{ is an irreducible bridge of span } L\},$$

$$\text{ir}^p(L; M, N) \stackrel{\text{def}}{=} \{\mathcal{S}^{(p)}[M, N] \text{ is an irreducible bridge of span } L\}$$

and

$$\lambda^p_{h,\beta}(L, N) \stackrel{\text{def}}{=} E_h[\exp(-\Phi^{(p)}_\beta(N)); \text{ir}^p(L; N)],$$

$$\Lambda^p_{h,\beta}(L) \stackrel{\text{def}}{=} \sum_{N \in \mathbb{N}^p} \lambda^p_{h,\beta}(L, N).$$

REMARK. In contrast to the definition of a bridge, the definition of the irreducibility of a bridge cannot be reduced to subpaths of the single random walks $\mathcal{S}^1, \ldots, \mathcal{S}^p$. This is manifested by the fact that the analogue of equation (2.13), which is a statement for bridges, becomes an inequality for irreducible bridges: for all $L \in \mathbb{N}$, we have

$$\overline{\Lambda}^p_{h,\beta}(L) \geq (\overline{\Lambda}_{h,\beta}(L))^p,$$

where $\overline{\Lambda}_{h,\beta}(L) \stackrel{\text{def}}{=} \overline{\Lambda}^1_{h,\beta}(L)$.

The following proposition states the so-called *renewal equation*, which provides access to the tools of renewal theory.

PROPOSITION 2.14. *For all $h, \beta \geq 0$ and $L \in \mathbb{N}$, we have*

(2.16) $$\overline{B}^p_{h,\beta}(L) = \sum_{k=1}^{L} \overline{\Lambda}^p_{h,\beta}(k) \overline{B}^p_{h,\beta}(L-k).$$

PROOF. For $\omega \in \text{br}^p(L; N)$, let $k \in \{1, \ldots, L\}$ denote the smallest breaking point for $\mathcal{S}^{(p)}[N](\omega)$. Then, there is a unique time $n \in \mathbb{N}^p$ with $n \leq N$ such that $\mathcal{S}^{(p)}[n](\omega)$ is an irreducible bridge of span $k$ and $\mathcal{S}^{(p)}[n, N](\omega)$ is a bridge of span $L - k$. We thus obtain

$$\text{br}^p(L; N) = \bigcup_{k=1}^{L} \bigcup_{n \in \mathbb{N}^p : n \leq N} \text{ir}^p(k; n) \cap \text{br}^p(L-k; n, N),$$

where the union is of disjoint sets. For $\omega \in \text{ir}^p(k; n) \cap \text{br}^p(L-k; n, N)$, we further have

$$\Phi^{(p)}_\beta(N, \omega) = \Phi^{(p)}_\beta(n, \omega) + \Phi^{(p)}_\beta(n, N, \omega),$$

by Lemma 2.1(b) applied to $\Phi^1_\beta, \ldots, \Phi^p_\beta$. By using the Markov property to renew the process $\mathcal{S}^{(p)}$ at time $n$, we thus have $\overline{B}^p_{h,\beta}(L)$ equal to

$$\sum_{k=1}^{L} \sum_{N \in \mathbb{N}^p} \sum_{n \in \mathbb{N}^p : n \leq N} E_h[e^{-\Phi^{(p)}_\beta(n)} 1_{\text{ir}^p(k;n)} e^{-\Phi^{(p)}_\beta(n,N)} 1_{\text{br}^p(L-k;n,N)}]$$

$$= \sum_{k=1}^{L} \sum_{N \in \mathbb{N}^p} \sum_{n \in \mathbb{N}^p : n \leq N} E_h[e^{-\Phi^{(p)}_\beta(n)} 1_{\text{ir}^p(k;n)}] E_h[e^{-\Phi^{(p)}_\beta(N-n)} 1_{\text{br}^p(L-k;N-n)}]$$

$$= \sum_{k=1}^{L} \overline{\Lambda}^p_{h,\beta}(k) \overline{B}^p_{h,\beta}(L-k). \qquad \square$$

REMARK. The subadditivity property of the sequence $(-\log \overline{B}(L))_{L \in \mathbb{N}}$, shown in Section 2.1 in a "straightforward" way, is also a consequence of the renewal equation for $p = 1$.

For $h, \beta \geq 0$, $L \in \mathbb{N}_0$ and $k \in \mathbb{N}$, we set

$$a_{h,\beta}^p(L) \stackrel{\text{def}}{=} \overline{B}_{h,\beta}^p(L) e^{\overline{m}_B^p(h,\beta)L} \quad \text{and} \quad \pi_{h,\beta}^p(k) \stackrel{\text{def}}{=} \overline{\Lambda}_{h,\beta}^p(k) e^{\overline{m}_B^p(h,\beta)k}.$$

LEMMA 2.15. *For any $h > 0$ and $\beta \geq 0$, we have*

(2.17) $$\sum_{k \in \mathbb{N}} \pi_{h,\beta}^p(k) = 1 \quad \text{and} \quad \sum_{k \in \mathbb{N}} k \pi_{h,\beta}^p(k) < \infty.$$

PROOF. For $L \in \mathbb{N}_0$, $k \in \mathbb{N}$ and $s \in [0, 1]$, set

$$A(s) \stackrel{\text{def}}{=} \sum_{L=0}^{\infty} a_{h,\beta}^p(L) s^L \quad \text{and} \quad P(s) \stackrel{\text{def}}{=} \sum_{k=1}^{\infty} \pi_{h,\beta}^p(k) s^k.$$

The renewal equation implies

$$A(s) = 1 + \sum_{L=1}^{\infty} a_{h,\beta}^p(L) s^L$$

$$= 1 + \sum_{L=1}^{\infty} \sum_{k=1}^{L} \pi_{h,\beta}^p(k) s^k a_{h,\beta}^p(L-k) s^{L-k}$$

$$= 1 + P(s) A(s).$$

By (2.14), we have $a_{h,\beta}^p(L) \leq 1$ for all $L \in \mathbb{N}_0$ and therefore $A(s) < \infty$ for $s \in [0, 1)$. As a consequence,

$$P(1) = \lim_{s \uparrow 1} P(s) = \lim_{s \uparrow 1} \frac{A(s) - 1}{A(s)} \leq 1.$$

We have thus shown that $\pi_{h,\beta}^p(k)$ is a (nonperiodic) subprobability sequence with renewal sequence $a_{h,\beta}^p(L)$. The first equation in (2.17) states that $\pi_{h,\beta}^p(k)$ is recurrent, which is equivalent to $A(1) = \infty$. If $\pi_{h,\beta}^p(k)$ is now recurrent, then the renewal theorem (see, e.g., Theorem 4.2.2 in [11]) yields

$$\lim_{L \to \infty} a_{h,\beta}^p(L) = \frac{1}{\sum_{k=1}^{\infty} k \pi_{h,\beta}^p(k)}.$$

The lemma thus follows from the estimate for $B_{h,\beta}^p(L)$ in (2.15), which states that $a_{h,\beta}^p(L) \geq \alpha(h)^p$ is true for all $L \in \mathbb{N}$. □

By (2.14) and Lemma 2.15, we know that $\overline{\Lambda}^p_{h,\beta}(L)$ decays to zero faster than $\overline{B}^p_{h,\beta}(L)$ when $L \to \infty$. In the next subsection, we will show that the difference in the decay velocity is even exponential. For this purpose, we will need the following result, stating that long intervals are unlikely to be free of breaking points. More precisely, we suppose that $\Delta, L \in \mathbb{N}$ with $L \geq 2\Delta$, $M \in \mathbb{N}_0^p$ and $N \in \mathbb{N}^p$ with $M < N$, and define

$$\mathrm{br}^{*p}_\Delta(L; M, N) \stackrel{\mathrm{def}}{=} \{\mathcal{S}^{(p)}[M, N] \text{ is a bridge of span } k, \text{ of which}$$
$$S^1(M^1) + \Delta, \ldots, S^1(N^1) - \Delta \text{ are no breaking points}\}$$

and $\mathrm{br}^{*p}_\Delta(L; N) \stackrel{\mathrm{def}}{=} \mathrm{br}^{*p}_\Delta(L; 0, N)$. For $h > 0$ and $\beta \geq 0$, we further set

$$b^{*p}_{\Delta,h,\beta}(L; N) \stackrel{\mathrm{def}}{=} \sum_{N \in \mathbb{N}^p} E_h[\exp(-\Phi^{(p)}_\beta(N)); \mathrm{br}^{*p}_\Delta(L; N)],$$

$$\overline{B}^{*p}_{\Delta,h,\beta}(L) \stackrel{\mathrm{def}}{=} \sum_{N \in \mathbb{N}^p} b^{*p}_{\Delta,h,\beta}(L; N).$$

LEMMA 2.16. *For any $h > 0$ and $\beta \geq 0$, there exists a decreasing function $\varepsilon^p_{h,\beta} : \mathbb{R}^+ \to \mathbb{R}^+$ such that $\lim_{T \to \infty} \varepsilon^p_{h,\beta}(T) = 0$ and*

$$\overline{B}^{*p}_{\Delta,h,\beta}(L) \leq \varepsilon^p_{h,\beta}(L - 2\Delta) e^{-\overline{m}^p_B(h,\beta)L}$$

*for all $\Delta, L \in \mathbb{N}$ with $L \geq 2\Delta$.*

PROOF. By the Markov property and Lemma 2.1(b), with $\ell$ denoting the largest breaking point smaller than $\Delta$ (or $\ell = 0$ if there is no such point) and $k$ denoting the smallest breaking point greater than $L - \Delta$, we have

$$\overline{B}^{*p}_{\Delta,h,\beta}(L) = \sum_{\ell=0}^{\Delta-1} \sum_{k=L-\Delta+1}^{L} \overline{B}^p_{h,\beta}(\ell) \overline{\Lambda}^p_{h,\beta}(k-\ell) \overline{B}^p_{h,\beta}(L-k)$$
$$= \sum_{\tilde{\ell}=1}^{\Delta} \sum_{\tilde{k}=T+1}^{L-\Delta} \overline{B}^p_{h,\beta}(\Delta - \tilde{\ell}) \overline{\Lambda}^p_{h,\beta}(\tilde{k} + \tilde{\ell}) \overline{B}^p_{h,\beta}(L - \tilde{k} - \Delta)$$

with $T = L - 2\Delta$, $\tilde{k} = k - \Delta$ and $\tilde{\ell} = \Delta - \ell$. Now, we set

(2.18) $$\varepsilon^p_{h,\beta}(T) \stackrel{\mathrm{def}}{=} \sum_{j=T+2}^{\infty} j \pi^p_{h,\beta}(j).$$

From Lemma 2.15, we know that $\lim_{T \to \infty} \varepsilon^p_{h,\beta}(T) = 0$. Moreover, we have

$$e^{\overline{m}^p_B(h,\beta)L} \overline{B}^{*p}_{\Delta,h,\beta}(L) \leq \sum_{\tilde{\ell}=1}^{\infty} \sum_{\tilde{k}=T+1}^{\infty} \pi^p_{h,\beta}(\tilde{k} + \tilde{\ell})$$

$$= \sum_{\tilde{\ell}=1}^{\infty} \sum_{j=T+\tilde{\ell}+1}^{\infty} \pi_{h,\beta}^p(j)$$

$$= \sum_{j=T+2}^{\infty} \sum_{\tilde{\ell}=1}^{j-T-1} \pi_{h,\beta}^p(j)$$

$$= \sum_{j=T+2}^{\infty} (j-T-1)\pi_{h,\beta}^p(j) < \varepsilon_{h,\beta}^p(T),$$

which completes the proof. $\square$

2.3. *Separation of the masses.* In this subsection, we investigate the decay velocity of weighted irreducible bridges. In particular, we will show that they decay exponentially faster than bridges without the irreducibility restriction.

LEMMA 2.17. *For all $h, \beta \geq 0$, the mass*

$$\overline{m}_\Lambda^p(h,\beta) \stackrel{\text{def}}{=} \lim_{L \to \infty} \frac{-\log \overline{\Lambda}_{h,\beta}^p(L)}{L}$$

*of $\overline{\Lambda}_{h,\beta}^p$ exists in $[\varphi_\beta(1), \infty)$, is continuous as function on $\mathbb{R}^+ \times \mathbb{R}^+$ and satisfies*

(2.19) $$\overline{\Lambda}_{h,\beta}^p(L) \leq \frac{1}{p} e^{2(\varphi_\beta(1)+\lambda_h)} e^{-\overline{m}_\Lambda^p(h,\beta)L}$$

*for all $L \in \mathbb{N}$, where $\lambda_h = \log E_0[\exp(h \cdot S_1(1))]$.*

PROOF. For $i \in \{1, \ldots, p\}$, let $E_i \stackrel{\text{def}}{=} (\delta_{i1}, \ldots, \delta_{ip}) \in \mathbb{N}_0^p$. Then, for every $N \in \mathbb{N}^p$, the union

$$\bigcup_{i=1}^{p} \bigcup_{N-2E_i > M \in \mathbb{N}^p} \mathrm{ir}^p(L_1; M) \cap \{S_1^i(M^i+1) = L_1+1\}$$

$$\cap \{S_1^i(M^i+2) = L_1\} \cap \mathrm{ir}^p(L_2; M+2E_i, N)$$

consists of disjoint sets and is a subset of $\mathrm{ir}^p(L_1+L_2; N)$. For any $\omega$ in that union, a double application of Lemma 2.1(b) to the potentials $\Phi^1, \ldots, \Phi^p$ yields

$$\Phi_\beta^{(p)}(N,\omega) \leq \Phi_\beta^{(p)}(M,\omega) + 2\varphi_\beta(1) + \Phi_\beta^{(p)}(M+2E_i, N, \omega)$$

for the corresponding $i \in \{1, \ldots, p\}$ and $M \in \mathbb{N}^p$. By splitting over all possible values of $i$ and $M$, and renewing the process $\mathcal{S}^{(p)}$ at the times $M$ and

$M + 2E_i$, we therefore have

$$\overline{\Lambda}^p_{h,\beta}(L_1 + L_2)$$
$$\geq \sum_{i=1}^{p} \sum_{N \in \mathbb{N}^p} \sum_{N-2E_i > M \in \mathbb{N}^p} E_h[e^{-\Phi^{(p)}_\beta(M)} 1_{\mathrm{ir}^p(L_1;M)} e^{-2\varphi_\beta(1)} 1_{\{S^i_1(M^i+1)=L_1+1\}}$$
$$\times 1_{\{S^i_1(M^i+2)=L_1\}} e^{-\Phi^{(p)}_\beta(M+2E_i,N)}$$
$$\times 1_{\mathrm{ir}^p(L_2;M+2E_i,N)}]$$
$$= \sum_{i=1}^{p} \sum_{N \in \mathbb{N}^p} \sum_{N-2E_i > M \in \mathbb{N}^p} \lambda^p_{h,\beta}(L;M) e^{-2\varphi_\beta(1)} P_h[S_1(1) = 1]$$
$$\times P_h[S_1(1) = -1] \lambda^p_{h,\beta}(L; N-M-2E_i)$$
$$= \frac{p}{e^{2(\varphi_\beta(1)+\lambda_h)}} \overline{\Lambda}^p_{h,\beta}(L_1) \overline{\Lambda}^p_{h,\beta}(L_2),$$

where, at the last step, the concrete definition of $P_h$ is used. The existence of $\overline{m}_\Lambda(h,\beta)$ in $[-\infty, \infty)$, as well as the estimate in (2.19), now follows from the subadditive limit lemma applied to

$$-\log\bigg(\frac{p}{e^{2(\varphi_\beta(1)+\lambda_h)}} \overline{\Lambda}^p_{h,\beta}(L)\bigg), \qquad L \in \mathbb{N}.$$

The lower bound $\varphi_\beta(1)$ for the mass goes back to Remark 2.3. Finally, the continuity of $\overline{m}_\Lambda$ is derived by the same arguments that were used in Proposition 2.7 to show the continuity of $\overline{m}_B$. □

The main result of this section is the derivation of a gap between the exponential decay rates of bridges and irreducible bridges.

THEOREM 2.18. *For any $h > 0$ and $\beta \geq 0$, we have*
$$\overline{m}^p_\Lambda(h,\beta) > \overline{m}^p_B(h,\beta).$$

The strategy for the proof of Theorem 2.18 was introduced in [4] (or see the more polished version of it in Chapter 4 of [11]) in the case of a single random walk in absence of a potential. It was then extended in [14] to single random walks evolving under the "trap" potential

$$\Phi^{\mathrm{trap}}_\beta(N) \stackrel{\mathrm{def}}{=} \beta \sharp R(N).$$

Before we present the strategy for the proof, we introduce the essential concept of backtracks of bridges.

DEFINITION 2.19. Suppose that $L \in \mathbb{N}$ and $\omega \in \text{br}^p(L; M, N)$ for $M \in \mathbb{N}_0^p$ and $N \in \mathbb{N}^p$ with $M < N$. Assume, furthermore, that $j \in \{1, \ldots, p\}$ and $M^j \leq m < n < N^j$ for $m \in \mathbb{N}_0$ and $n \in \mathbb{N}$.

The subpath $\mathcal{S}^j[m, n](\omega)$ of $\mathcal{S}^j[M^j, N^j](\omega)$ is called a *j-backtrack* (or simply *backtrack*) of the bridge $\mathcal{S}^{(p)}[M, N](\omega)$ if

(i) $S_1^j(\mu_1, \omega) \leq S_1^j(m, \omega)$ for $\mu_1 = M^j + 1, \ldots, m$;
(ii) $S_1^j(n, \omega) \leq S_1^j(\nu, \omega) < S_1^j(m, \omega)$ for $\nu = m + 1, \ldots, n$;
(iii) $S_1^j(n, \omega) < S_1^j(\mu_2, \omega)$ for $\mu_2 = n + 1, \ldots, N^j$.

If this is the case, then the *span* of the backtrack $\mathcal{S}^j[m,n](\omega)$ is given by $S_1^j(m, \omega) - S_1^j(n, \omega)$. A backtrack $\mathcal{S}^j[m, n](\omega)$ is said to *cover* an integer $k$ if $S_1^j(n, \omega) \leq k < S_1^j(m, \omega)$ is valid.

REMARKS. (a) Condition (ii) says that a backtrack is itself a bridge in $\mathbb{Z}^d$, except that it goes "right-to-left" instead of "left-to-right." Conditions (i) and (iii) are maximality conditions. For two different $j$-backtracks, $\mathcal{S}^j[m_1, n_1](\omega)$ and $\mathcal{S}^j[m_2, n_2](\omega)$, they guarantee that the time intervals $\{m_1, \ldots, n_1\}$ and $\{m_2, \ldots, n_2\}$ do not intersect and that $n_1 < m_2$ is equivalent to $S_1^j(n_1, \omega) < S_1^j(n_2, \omega)$.

(b) A bridge $\mathcal{S}^{(p)}[M, N](\omega)$ is irreducible if and only if every integer $k$ with $S^1(M^1, \omega) < k < S^1(N^1, \omega)$ is covered by a backtrack of $\mathcal{S}^{(p)}[M, N](\omega)$. An integer $k$ may of course be covered by several backtracks.

We now present the strategy for the proof of Theorem 2.18. The aim is to find an upper bound for $\overline{\Lambda}_{h,\beta}^p(L)$. To this end, we fix a large integer $Q$ and split the interval $[0, L]$ for $L \gg Q$ into blocks (subintervals) of size $Q$. For an irreducible bridge of span $L$, we then look at the backtracks that cover the endpoints of these blocks and distinguish between the two following situations.

In the first situation, many of these endpoints are covered by only small backtracks (of span not larger than $\Lambda \ll Q$). Between such endpoints, the path consists of subbridges with large intervals being free of breaking points. This will allow a multiple application of Lemma 2.16.

In the other situation, some of the endpoints are covered by large backtracks. In that situation, the random walk must go "backward" often, which it does with small probability because of the drift. It is going to be important that the "total span" of these backtracks remains large enough with respect to the reduced number of points considered (i.e., the number of endpoints).

More precisely, we proceed as follows. Let $T$ and $\Delta$ be positive integers (to be specified) and set $Q \stackrel{\text{def}}{=} 2\Delta + T$. For large $L \in \mathbb{N}$, let $k = k(L)$ be the greatest integer less than or equal to $\frac{L}{Q} - 1$ and set

$$A \stackrel{\text{def}}{=} \{Q, 2Q, \ldots, kQ\}.$$

Now, let $B = \{b_1, \ldots, b_\tau\}$ with $b_1 < \cdots < b_\tau$ be a subset of $A$ and observe that
$$b_i - b_{i-1} \geq Q \quad \text{for } i = 1, \ldots, \tau + 1,$$
where $b_0 \stackrel{\text{def}}{=} 0$ and $b_{\tau+1} \stackrel{\text{def}}{=} L$. We introduce two further items of notation for particular $N$-step bridges of span $L$:

- for $\Delta > 0$, let
$$\text{ir}^p_{\Delta, B}(L; N) \in \mathcal{F}^{\otimes p}$$
denote the set of all $\omega \in \Omega^p$ for which $\mathcal{S}^{(p)}[N](\omega)$ is an irreducible bridge of span $L$ such that no point of $B$ is covered by a backtrack of $\mathcal{S}^{(p)}[N](\omega)$ with a span larger than $\Delta$;

- for $\sigma \in \mathbb{N}^B$, any pairwise disjoint decomposition $B^1, \ldots, B^p$ of $B$ (possibly with some of the $B^j$ being empty) and each $j \in \{1, \ldots, p\}$, let
$$\text{br}^{p,j}_{B^j, \sigma|_{B^j}}(L; N^j) \in \mathcal{F}^{\otimes p}$$
denote the set of all $\omega \in \Omega^p$ for which $\mathcal{S}^j[N^j](\omega)$ is a bridge of span $L$ such that each $b \in B^j$ is covered by a backtrack of $\mathcal{S}^j[N^j](\omega)$ with span $\sigma(b)$ which is not covering any other $a \in B^j \setminus \{b\}$, and let
$$\text{br}^p_{B^1, \ldots, B^p, \sigma}(L; N) \stackrel{\text{def}}{=} \bigcap_{j=1}^p \text{br}^{p,j}_{B^j, \sigma|_{B^j}}(L; N^j).$$

The following lemma realizes the aforementioned distinction on how the points of $A = \{Q, 2Q, \ldots, kQ\}$ are covered by backtracks.

LEMMA 2.20. *The event $\text{ir}^p(L; N)$ is contained in the union of*
$$\bigcup_{B \subset A : \sharp B \geq k/2} \text{ir}^p_{\Delta, B}(L; N)$$
*and*
$$\bigcup_{\substack{B \subset A : \sharp B \geq 1}} \bigcup_{\substack{\sigma \in \mathbb{N}^B : \\ \sum_{b \in B} \sigma(b) > k\Delta/2}} \bigcup_{\substack{B^1, \ldots, B^p \subset B \\ \text{pairwise disjoint} \\ \text{decomposition}}} \text{br}^p_{B^1, \ldots, B^p, \sigma}(L; N).$$

PROOF. Suppose that $\omega$ is in $\text{ir}^p(L; N)$, but not in
$$\bigcup_{B \subset A, \sharp B \geq k/2} \text{ir}^p_{\Delta, B}(L; N).$$
There then exists a collection of backtracks of $\mathcal{S}^{(p)}[N](\omega)$, each of them of a span greater than $\Delta$, which cover at least $k/2$ of the points in $A$. Although some of them may cover several points in $A$, the sum of their spans is still greater than $\Delta k/2$, since the distance between two points is at least $Q = T + 2\Delta$.

We now inductively construct the sets $B^1, \ldots, B^p$. For $j \in \{1, \ldots, p\}$, assume that $B^1, \ldots, B^{j-1}$ are already constructed. For each $j$-backtrack from our collection, we then add to $B^j$ a single point from the complement of $B^1 \cup \cdots \cup B^{j-1}$ which is covered by this $j$-backtrack, but not covered by any other $j$-backtrack (regardless of whether it is covered by $i$-backtracks for $i \neq j$). If there is no such point, we remove this particular $j$-backtrack from the collection. The remaining backtracks still cover the same points in $A$, so the sum of their spans is still larger than $\Delta k/2$.

By this construction, the sets $B^1, \ldots, B^p$ are pairwise disjoint. Moreover, for any $j \in \{1, \ldots, p\}$, the set $B^j$ has the property, that each of its points is covered by exactly one of the remaining $j$-backtracks. Consequently, if we set $B \stackrel{\text{def}}{=} B^1 \cup \cdots \cup B^p$, then there is a $\sigma \in \mathbb{N}^B$ with $\sum_{b \in B} \sigma(b) > \Delta k/2$ and $\omega \in \operatorname{br}^p_{B^1,\ldots,B^p,\sigma}(L; N)$. □

For $B = \{b_1, \ldots, b_\tau\} \subset A$ with $b_1 < \cdots < b_\tau$, we set

$$\overline{\Lambda}^p_{h,\beta}[\operatorname{ir}^p_{\Delta,B}(L;N)] \stackrel{\text{def}}{=} \sum_{N \in \mathbb{N}^p} E_h[\exp(-\Phi^{(p)}_\beta(N)); \operatorname{ir}^p_{\Delta,B}(L;N)].$$

LEMMA 2.21. *For any $h > 0$ and $\beta \geq 0$, we have*

$$\overline{\Lambda}^p_{h,\beta}[\operatorname{ir}^p_{\Delta,B}(L;N)] \leq e^{-\overline{m}^p_B(h,\beta)L}(\varepsilon^p_{h,\beta}(T)\alpha(h)^{-p})^{\tau+1},$$

*where $\alpha(h)$ is defined in (2.5) and $\varepsilon^p_{h,\beta}(T)$ is defined in (2.18).*

PROOF. Suppose that $\omega \in \operatorname{ir}^p_{\Delta,B}(L;N)$ and recall that $b_0 = 0$ and $b_{\tau+1} = L$. For each $i \in \{1, \ldots, \tau+1\}$, let $m_{i-1} = (m^1_{i-1}, \ldots, m^p_{i-1}) \in \mathbb{N}_0^p$ and $n_i = (n^1_i, \ldots, n^p_i) \in \mathbb{N}^p$ be given by

$$m^j_{i-1} \stackrel{\text{def}}{=} \min\{\mu \in \{0, \ldots, N^j\} : S^j_1(\mu', \omega) > b_{i-1} \text{ for } \mu < \mu' \leq N^j\},$$

$$n^j_i \stackrel{\text{def}}{=} \max\{\nu \in \{1, \ldots, N^j\} : S^j_1(\nu', \omega) \leq b_i \text{ for } 0 < \nu' \leq \nu\}$$

for $j = 1, \ldots, p$. It is also convenient to choose $m_{\tau+1} \stackrel{\text{def}}{=} L$.

Since $b_1, \ldots, b_\tau$ are not covered by backtracks of a span greater than $\Delta$ and because the distance between these points is at least $Q = 2\Delta + T$, it is clear that

$$S^j_1(\mu, \omega) \leq b_{i-1} + \Delta < b_i - \Delta < S^j_1(\nu, \omega)$$

for $i = 1, \ldots, \tau + 1$, $j = 1, \ldots, p$ and $1 \leq \mu \leq m_{i-1}^j < n_i^j < \nu \leq N^j$. Therefore, since $\mathcal{S}^{(p)}[N](\omega)$ is irreducible, we know that the sub-bridges $\mathcal{S}^{(p)}[m_{i-1}, n_i](\omega)$ already contain backtracks that cover the points $b_{i-1} + \Delta, \ldots, b_i - \Delta$. That is, we have

$$\omega \in \bigcap_{i=1}^{\tau+1} \mathrm{br}_\Delta^{*p}(b_i - b_{i-1}; m_{i-1}, n_i),$$

where $\mathrm{br}_\Delta^{*p}$ was introduced at the end of Section 2.2. An inductive application of Lemma 2.1(c) to the potentials $\Phi^1, \ldots, \Phi^p$ further yields

$$\Phi_\beta^{(p)}(N, \omega) \geq \sum_{i=1}^{\tau+1} \Phi_\beta^{(p)}(m_{i-1}, n_i, \omega).$$

As a consequence, by the Markov property and Lemma 2.16, an upper bound for $\overline{\Lambda}_{h,\beta}^p[\mathrm{ir}_{\Delta,B}^p(L; N)]$ is given by

$$\sum_{N \in \mathbb{N}^p} \sum_{\substack{n_1, m_1, n_2, \ldots, n_\tau, m_\tau \in \mathbb{N}^p: \\ 0 = m_0 < n_1 < \cdots < m_\tau < n_{\tau+1} = m_{\tau+1} = N}} E_h \left[ \prod_{i=1}^{\tau+1} e^{-\Phi_\beta^{(p)}(m_{i-1}, n_i)} 1_{\mathrm{br}_\Delta^{*p}(b_i - b_{i-1}; m_{i-1}, n_i)} \right.$$

$$\left. \times 1_{\{S_1^j(m_i^j) = S_1^j(n_i^j) \text{ for } j=1, \ldots, p\}} \right]$$

$$= \sum_{N \in \mathbb{N}^p} \sum_{\substack{n_1, m_1, n_2, \ldots, n_\tau, m_\tau \in \mathbb{N}^p \\ 0 = m_0 < n_1 < \cdots < m_\tau < n_{\tau+1} = m_{\tau+1} = N}} \prod_{i=1}^{\tau+1} b_{\Delta,h,\beta}^{*p}(b_i - b_{i-1}; n_i - m_{i-1})$$

$$\times P_h[S_1^j(m_i^j - n_i^j) = 0 \text{ for } j = 1, \ldots, p]$$

$$= \left( \prod_{i=1}^{\tau+1} \overline{B}_{\Delta,h,\beta}^{*p}(b_i - b_{i-1}) \right) \left( \sum_{m \in \mathbb{N}^p} \prod_{j=1}^p P_h[S_1^j(m^j) = 0] \right)^\tau$$

$$\leq \left( \prod_{i=1}^{\tau+1} \varepsilon_{h,\beta}^p(b_i - b_{i-1} - 2\Delta) \right) e^{-\overline{m}_B^p(h,\beta)L} \frac{1}{\alpha(h)^{p\tau}},$$

where, at the last step, (2.6) is used to identify $\alpha(h)$. The lemma now follows from the facts that $\alpha(h) \leq 1$ and $b_i - b_{i-1} \geq Q = T + 2\Delta$ for $i = 1, \ldots, \tau + 1$. □

For $B = \{b_1, \ldots, b_\tau\} \subset A$ with $b_1 < \cdots < b_\tau$, and a pairwise disjoint decomposition $B^1, \ldots, B^p$ of $B$, we define

$$\overline{\Lambda}_{h,\beta}^p[\mathrm{br}_{B^1, \ldots, B^p, \sigma}^p(L; N)] \stackrel{\text{def}}{=} \sum_{N \in \mathbb{N}^p} E_h[\exp(-\Phi_\beta^{(p)}(N)); \mathrm{br}_{B^1, \ldots, B^p, \sigma}^p(L; N)].$$

LEMMA 2.22. *For any $h > 0$, $\beta \geq 0$ and $\sigma \in \mathbb{N}^B$, we have*

$$\sum_{\substack{B^1,\ldots,B^p \subset B \\ \text{pairwise disjoint} \\ \text{decomposition}}} \overline{\Lambda}^p_{h,\beta}[\text{br}^p_{B^1,\ldots,B^p,\sigma}(L;N)] \leq \psi^p(h)^\tau e^{-\overline{m}^p_B(h,\beta)L} e^{-h \sum_{b \in B} \sigma(b)},$$

*where $\alpha(h)$ is defined in (2.5) and $\psi^p(h) \stackrel{\text{def}}{=} p/\alpha(h)^2 h$.*

PROOF. We first restrict to the case $p = 1$. For $L \in \mathbb{N}_0$ and $M, N \in \mathbb{N}_0$ with $M \leq N$, we define

$$\text{br}'(L;N) \stackrel{\text{def}}{=} \{-\mathcal{S}[N] \text{ is a bridge of span } L\},$$

$$\text{br}'(L;M,N) \stackrel{\text{def}}{=} \{-\mathcal{S}[M,N] \text{ is a bridge of span } L\}.$$

By the Markov property of $\mathcal{S}$, it is easily seen that

$$P_h[\text{br}'(L;N)] = P_h[-L = S_1(N) < S_1(n) \leq 0 \text{ for } 0 \leq n < N].$$

By Remark 2.2, we thus have

$$(2.20) \qquad \overline{B}'_{h,0}(L) \stackrel{\text{def}}{=} \sum_{N=1}^\infty P_h[\text{br}'(L;N)] \leq P_h[H_{-L} < \infty] = e^{-2hL},$$

where $H_{-L} = \inf\{n \in \mathbb{N} : S_1(n) = -L\}$.

Suppose now that $\omega \in \text{br}^1_{B,\sigma}(L;N)$ and recall that $b_0 = 0$ and $b_{\tau+1} = L$. For each $i \in \{1,\ldots,\tau\}$, by the definition of $\text{br}^1_{B,\sigma}(L;N)$, there is at least one backtrack $\mathcal{S}[s_i, t_i](\omega)$ of $\mathcal{S}[N](\omega)$ of span $\sigma(b_i)$ satisfying

$$(2.21) \qquad b_{i-1} < S_1(t_i,\omega) \leq b_i < S_1(s_i,\omega) \leq b_{i+1}.$$

It is also convenient to choose $s_{\tau+1} \stackrel{\text{def}}{=} t_{\tau+1} \stackrel{\text{def}}{=} L$ and $\sigma_{\tau+1} \stackrel{\text{def}}{=} 0$, in order to have $P_h[\text{br}'(\sigma_{\tau+1}; s_{\tau+1}, t_{\tau+1})] = 0$.

For $i \in \{1,\ldots,\tau+1\}$, we now define

$$m_{i-1} \stackrel{\text{def}}{=} \min\{\mu \in \{0,\ldots,N\} : S_1(\mu',\omega) > b_{i-1} \text{ for } \mu < \mu' \leq N\},$$

$$n_i \stackrel{\text{def}}{=} \max\{\nu \in \{m_{i-1}+1,\ldots,N\} : S_1(\nu',\omega) \leq b_i \text{ for } \mu < \nu' \leq \nu\},$$

such that $\mathcal{S}[m_{i-1}, n_i](\omega)$ is a bridge of span $b_i - b_{i-1}$. By (2.21) and the maximality conditions for backtracks, we have $t_i < s_{i+1}$ and thus

$$t_i \leq m_i < n_{i+1} \leq s_{i+1}$$

for $i = 1,\ldots,\tau$, by the definition of $n_{i+1}$ and (2.21) again. Moreover, since

$$\Phi_\beta(N,\omega) \geq \sum_{i=1}^{\tau+1} \Phi_\beta(m_{i-1}, n_i, \omega)$$

by Lemma 2.1(c), the Markov property of the random walk yields

$$\overline{\Lambda}^1_{h,\beta}[\text{br}^1_{B,\sigma}(L;N)]$$

$$\leq \sum_{N\in\mathbb{N}} \sum_{\substack{0=m_0<n_1<s_1<t_1\leq m_1<\cdots \\ \cdots<n_{\tau+1}=s_{\tau+1}=t_{\tau+1}=m_{\tau+1}=N}} \sum_{\substack{q_1,\ldots,q_\tau\in\mathbb{N}, q_{\tau+1}=0, \tilde{q}_1,\ldots,\tilde{q}_\tau\in\mathbb{N}_0: \\ q_1+\tilde{q}_1=\sigma(b_1),\ldots,q_\tau+\tilde{q}_\tau=\sigma(b_\tau)}}$$

$$E_h\left[\prod_{i=1}^{\tau+1} e^{-\Phi_\beta(m_{i-1},n_i)} 1_{\text{br}(b_i-b_{i-1};m_{i-1},n_i)} 1_{\{S_1(s_i)=S_1(n_i)+q_i\}}\right.$$

$$\left.\times 1_{\text{br}'(\sigma(b_i);s_i,t_i)} 1_{\{S_1(m_i)=S_1(t_i)+\tilde{q}_i\}}\right]$$

$$= \sum_{N\in\mathbb{N}} \sum_{0=m_0<n_1<s_1<t_1\leq m_1<\cdots<n_{\tau+1}=N}$$

$$\left(\prod_{i=1}^{\tau+1} b_{h,\beta}(L;n_i-m_{i-1})\right)\left(\prod_{i=1}^{\tau} P_h[\text{br}'(\sigma(b_i);t_i-s_i)]\right)$$

$$\times \left(\prod_{i=1}^{\tau} \sum_{q\in\mathbb{N},\tilde{q}\in\mathbb{N}_0:q+\tilde{q}=\sigma(b_i)} P_h[S_1(s_i-n_i)=q]P_h[S_1(m_i-t_i)=\tilde{q}]\right)$$

$$\leq \left(\prod_{i=1}^{\tau+1} \overline{B}_{h,\beta}(L)\right)\left(\prod_{i=1}^{\tau} \overline{B}'_{h,0}(L)\right)$$

$$\times \left(\prod_{i=1}^{\tau} \sum_{q\in\mathbb{N},\tilde{q}\in\mathbb{N}_0:q+\tilde{q}=\sigma(b_i)} \sum_{k=0}^{\infty} P_h[S_1(k)=q] \sum_{\tilde{k}=1}^{\infty} P_h[S_1(\tilde{k})=\tilde{q}]\right)$$

$$= \alpha(h)^{-2\tau}\left(\prod_{i=1}^{\tau+1} \overline{B}_{h,\beta}(b_i-b_{i-1})\right)\left(\prod_{i=1}^{\tau} \overline{B}'_{h,0}(\sigma(b_i))\right)\left(\prod_{i=1}^{\tau} \sigma(b_i)\right),$$

where, in the last step, equation (2.6) is used to identify $\alpha(h)$. From Proposition 2.7 and (2.20), we then obtain

$$\overline{\Lambda}^1_{h,\beta}[\text{br}^1_{B,\sigma}(L;N)] \leq \alpha(h)^{-2\tau} e^{-\overline{m}_B(h,\beta)L} \prod_{i=1}^{\tau} \sigma(b_i) e^{-2h\sigma(b_i)}$$

$$\leq (\alpha(h)^2 h)^{-\tau} e^{-\overline{m}_B(h,\beta)L} e^{-h\sum_{i=1}^{\tau}\sigma(b_i)},$$

where, in the second estimate, the elementary inequality $xe^{-x}\leq 1$ for $x\geq 0$ is used.

We now proceed to the case of arbitrary $p\in\mathbb{N}$. We consider a pairwise disjoint decomposition $B^1,\ldots,B^p$ of $B=\{b_1,\ldots,b_\tau\}$. By the independence of the walks $\mathcal{S}^1,\ldots,\mathcal{S}^p$, we then have

$$\overline{\Lambda}^p_{h,\beta}[\mathrm{br}^p_{B^1,\ldots,B^p,\sigma}(L;N)] = \sum_{N\in\mathbb{N}^p}\prod_{j=1}^p E_h[\exp(-\Phi^j(N^j));\mathrm{br}^{p,j}_{B^j,\sigma|_{B^j}}(L;N)]$$

$$= \prod_{j=1}^p \overline{\Lambda}^1_{h,\beta}[\mathrm{br}^1_{B^j,\sigma|_{B^j}}(L;N)]$$

$$\leq (\alpha(h)^2 h)^\tau e^{-p\overline{m}_B(h,\beta)L} e^{-h\sum_{b\in B}\sigma(b)},$$

where, in the last step, the estimate for a single random walk is used. Since there are $p^\tau$ pairwise disjoint decompositions $B^1,\ldots,B^p$ of $B$, the lemma now follows from $\overline{m}^p_B(h,\beta) = p\overline{m}_B(h,\beta)$ in Remark 2.12. $\square$

PROOF OF THEOREM 2.18. By means of Lemma 2.20 and Lemma 2.22, we finally have the necessary tools to prove the mass gap for irreducible bridges. We first fix $T\in\mathbb{N}$ and $\Delta\in\mathbb{N}$ large enough such that we have

$$2(\alpha(h)^{-p}\varepsilon^p_{h,\beta}(T))^{1/2} \leq \frac{1}{2} \quad\text{and}\quad e^{-\Delta h/4}\left(1+\frac{\psi^p(h)}{1-e^{-h/2}}\right) \leq \frac{1}{2},$$

and we set $Q \stackrel{\mathrm{def}}{=} 2\Delta + T$. Since there are $2^k$ subsets of $A = \{Q,\ldots,kQ\}$, we have

$$\sum_{B\subset A:\sharp B\geq k/2}\overline{\Lambda}^p_{h,\beta}[\mathrm{ir}^p_{\Delta,B}(L;N)] \leq e^{-\overline{m}^p_B(h,\beta)L} 2^k(\alpha(h)^{-p}\varepsilon^p_{h,\beta}(T))^{1+k/2}$$

$$< e^{-\overline{m}^p_B(h,\beta)L} 2^{-k}$$

by Lemma 2.21. Moreover, since there are $\binom{k}{\tau}$ subsets $B\subset A$ with $\sharp B = \tau$, we obtain from Lemma 2.22 that

$$\sum_{B\subset A:\sharp B\geq 1}\sum_{\substack{\sigma\in\mathbb{N}^B:\\ \sum_{b\in B}\sigma(b)>k\Delta/2}}\sum_{\substack{B^1,\ldots,B^p\subset B\\ \text{pairwise disjoint}\\ \text{decomposition}}} \overline{\Lambda}^p_{h,\beta}[\mathrm{br}^p_{B^1,\ldots,B^p,\sigma}(L;N)]$$

$$\leq e^{-\overline{m}^p_B(h,\beta)L}\sum_{\tau=1}^k\binom{k}{\tau}\psi^p(h)^\tau \sum_{\substack{\sigma_1,\ldots,\sigma_\tau\in\mathbb{N}:\\ \sigma_1+\cdots+\sigma_\tau>k\Delta/2}} e^{-h(\sigma_1+\cdots+\sigma_\tau)}$$

$$\leq e^{-\overline{m}^p_B(h,\beta)L}e^{-k\Delta h/4}\sum_{\tau=1}^k\binom{k}{\tau}\psi^p(h)^\tau \sum_{\sigma_1,\ldots,\sigma_\tau\in\mathbb{N}} e^{(-h/2)(\sigma_1+\cdots+\sigma_\tau)}$$

$$\leq e^{-\overline{m}^p_B(h,\beta)L}e^{-k\Delta h/4}\sum_{\tau=1}^k\binom{k}{\tau}\left(\frac{\psi^p(h)}{1-e^{-h/2}}\right)^\tau$$

$$\leq e^{-\overline{m}_B^p(h,\beta)L} e^{-k\Delta h/4} \left(1 + \frac{\psi^p(h)}{1-e^{-h/2}}\right)^k$$

$$\leq e^{-\overline{m}_B^p(h,\beta)L} 2^{-k}.$$

Finally, we apply Lemma 2.20 and the estimate $k \geq \frac{L}{Q} - 2$ to obtain

$$\overline{\Lambda}_{h,\beta}^p(L) \leq 2e^{-\overline{m}_B^p(h,\beta)L} 2^{-k} \leq 8e^{-\overline{m}_B^p(h,\beta)L} 2^{-L/Q}$$

and therefore

$$\overline{m}_\Lambda^p(h,\beta) = \lim_{L\to\infty} \frac{-\log \overline{\Lambda}_{h,\beta}^p(L)}{L} \geq \overline{m}_B^p(h,\beta) + \frac{\log 2}{Q} > \overline{m}_B^p(h,\beta),$$

which proves the theorem. □

**3. Fixed number of steps.** In this section, we consider finite random walks $\mathcal{S}[N]$ for fixed $N \in \mathbb{N}$, again evolving under the influence of the path potential

$$\Phi_\beta(N) = \sum_{x \in \mathbb{Z}^d} \varphi_\beta(\ell_x(N)), \qquad N \in \mathbb{N}.$$

As introduced at the beginning of Section 2, we assume $\varphi_\beta$ to be given by

$$\varphi_\beta(t) = \varphi(\beta t), \qquad t \in \mathbb{R}^+,$$

where $\varphi : \mathbb{R}^+ \to \mathbb{R}^+$ is a concave increasing function satisfying

$$\lim_{t\to 0} \varphi(t) = \varphi(0) = 0,$$

as well as $\lim_{t\to\infty} \varphi(t) = \infty$ and $\lim_{t\to\infty} \varphi(t)/t = 0$.

3.1. *Masses for paths and bridges.* In the present setting, for $h, \beta \geq 0$, the generalization of the annealed partition function $Z_{h,\beta}^{\mathrm{an}}$ from Section 1 is given by

$$G_{h,\beta}(N) \stackrel{\mathrm{def}}{=} E_h[\exp(-\Phi_\beta(N))], \qquad N \in \mathbb{N}.$$

REMARK 3.1. For any $h, \beta \geq 0$ and $N \in \mathbb{N}$, we have

$$e^{-\varphi_\beta(1)N} \leq G_{h,\beta}(N) \leq 1,$$

where the lower estimate follows from Lemma 2.1(a).

We are interested in the limiting exponential behavior of $G_{h,\beta}(N)$. The existence of an associated mass was part of Theorem A (which, in the original paper [7], is shown in the present, more general setting). As we show next, it can also be obtained in a straightforward way by the subadditive limit lemma, which additionally delivers a bound for the speed of convergence.

PROPOSITION 3.2. *For any $h, \beta \geq 0$, the mass*

$$m_G(h, \beta) \stackrel{\text{def}}{=} \lim_{N \to \infty} \frac{-\log G_{h,\beta}(N)}{N}$$

*of $G_{h,\beta}$ exists in $[0, \varphi_\beta(1)]$ and is continuous as a function on $\mathbb{R}^+ \times \mathbb{R}^+$, and*

$$G_{h,\beta}(N) \leq e^{-m_G(h,\beta)N}$$

*is valid for all $N \in \mathbb{N}$.*

PROOF. By Lemma 2.1(b), for any $N_1, N_2 \in \mathbb{N}$, we have

$$\begin{aligned}
G_{h,\beta}(N_1 + N_2) &= E_h[e^{-\Phi_\beta(N_1+N_2)}] \\
&\geq E_h[e^{-\Phi_\beta(N_1)-\Phi_\beta(N_1,N_1+N_2)}] = G_{h,\beta}(N_1)G_{h,\beta}(N_2),
\end{aligned} \tag{3.1}$$

where, in the last step, the Markov property is used to renew the random walk at time $N_1$. Therefore, the existence of the mass $m_G$ and the estimate in (3.2) are consequences of the subadditive limit lemma applied to the sequence $(-\log G_{h,\beta}(N))_{N \in \mathbb{N}}$ and the bounds for the mass follow from Remark 3.1. Finally, the continuity of $m_G$ is obtained by similar (but slightly simpler) arguments as used to prove the continuity of $\overline{m}_B$ in Proposition 2.7. □

As in the point-to-hyperplane setting of Section 2, it is convenient to introduce $N$-step bridges. For $M, N \in \mathbb{N}_0$ with $M \leq N$, we define

$$\text{Br}(N) \stackrel{\text{def}}{=} \{\mathcal{S}[N] \text{ is a bridge}\},$$

$$\text{Br}(M, N) \stackrel{\text{def}}{=} \{\mathcal{S}[M, N] \text{ is a bridge}\}.$$

For $h, \beta \geq 0$ and $N \in \mathbb{N}_0$, we further set

$$B_{h,\beta}(N) \stackrel{\text{def}}{=} E_h[\exp(-\Phi_\beta(N)); \text{Br}(N)].$$

Observe that we have

$$\text{Br}(N) = \bigcup_{L \in \mathbb{N}_0} \text{br}(L; N) \quad \text{and} \quad B_{h,\beta}(N) = \sum_{L \in \mathbb{N}_0} b_{h,\beta}(L; N),$$

where the union is of disjoint sets and where

$$\text{br}(L; N) = \{\mathcal{S}[N] \text{ is a bridge of span } L\},$$
$$b_{h,\beta}(L; N) = E_h[\exp(-\Phi_\beta(N)); \text{br}(L; N)]$$

were introduced in Section 2.1.

REMARK 3.3. For any $h, \beta \geq 0$ and $N \in \mathbb{N}$, we have
$$e^{-\varphi_\beta(1)N} B_{h,0}(N) \leq B_{h,\beta}(N) \leq G_{h,\beta}(N),$$
where the lower estimate follows from Lemma 2.1(a).

Our interest lies with the ballistic regime $\{(h,\beta) \in (0,\infty)^2 : \beta < \beta_c(h)\}$, where the critical parameter $\beta_c(h)$, as introduced in Remark 2.10, is determined by $\overline{m}(0, \beta_c(h)) = h$. In our first result on $N$-step bridges, however, it is not necessary to restrict to the ballistic phase.

LEMMA 3.4. *For any $h, \beta \geq 0$ and $N_1, N_2 \in \mathbb{N}$, we have*

(3.2) $$B_{h,\beta}(N_1 + N_2) \geq B_{h,\beta}(N_1) B_{h,\beta}(N_2),$$

(3.3) $$B_{h,\beta}(N_1 + N_2) \leq K_{h,\beta} B_{h,\beta}(N_1) B_{h,\beta}(N_2),$$

*where*

(3.4) $$K_{h,\beta} \stackrel{\text{def}}{=} \left( \sum_{n=0}^{\infty} \frac{P_h[S_1(n) = 0]}{B_{h,\beta}(n)} \right)^2.$$

PROOF. By the definition of a bridge, it is plain that
$$\text{Br}(N_1 + N_2) \supset \text{Br}(N_1) \cap \text{Br}(N_1, N_1 + N_2).$$
From Lemma 2.1(b), we thus obtain
$$B_{h,\beta}(N_1 + N_2) \geq E_h[e^{-\Phi_\beta(N_1)} 1_{\text{Br}(N_1)} e^{-\Phi_\beta(N_1, N_1+N_2)} 1_{\text{Br}(N_1, N_1+N_2)}]$$
$$= B_{h,\beta}(N_1) B_{h,\beta}(N_2),$$
where, in the second step, the Markov property is used to renew the random walk at time $N_1$.

For the upper estimate, observe that we have
$$\text{Br}(N_1 + N_2) \subset \bigcup_{m_1=1}^{N_1} \bigcup_{m_2=N_1}^{N_1+N_2} \text{Br}(m_1) \cap \{S_1(N_1) = S_1(m_1)\}$$
$$\cap \{S_1(m_2) = S_1(N_1)\} \cap \text{Br}(m_2, N_1 + N_2),$$
where $m_1$ and $m_2$ are the times of the first and last visits, respectively, of $S[N_1 + N_2]$ to the hyperplane $\mathcal{H}_{S_1(N_1)}$. By splitting over all possible values of $m_1$ and $m_2$, and by applying Lemma 2.1(c) and the Markov property to renew the walk at these times, we obtain
$$B_{h,\beta}(N_1 + N_2)$$
$$\leq \sum_{m_1=1}^{N_1} \sum_{m_2=N_1}^{N_1+N_2} E_h[e^{-\Phi_\beta(m_1)} 1_{\text{Br}(m_1)} 1_{\{S_1(N_1)=S_1(m_1)\}}$$

$$\times 1_{\{S_1(m_2)=S_1(N_1)\}} e^{-\Phi_\beta(m_2,N_1+N_2)} 1_{\mathrm{Br}(m_2,N_1+N_2)}]$$

$$= \sum_{m_1=1}^{N_1} \sum_{m_2=N_1}^{N_1+N_2} B_{h,\beta}(m_1) P_h[S_1(N_1-m_1)=0]$$

$$\times P_h[S_1(m_2-N_1)=0] B_{h,\beta}(N_1+N_2-m_2)$$

$$= B_{h,\beta}(N_1) B_{h,\beta}(N_2) \left( \sum_{n_1=0}^{N_1-1} \frac{P_h[S_1(n_1)=0]}{B_{h,\beta}(n_1)} \right) \left( \sum_{n_2=0}^{N_2} \frac{P_h[S_1(n_2)=0]}{B_{h,\beta}(n_2)} \right),$$

where $n_1$ replaces $N_1 - m_1$ and $n_2$ replaces $m_2 - N_1$, and where, in the last step, (3.2) is used. □

PROPOSITION 3.5. *For any $h, \beta \geq 0$, the mass*

$$m_B(h,\beta) \stackrel{\mathrm{def}}{=} \lim_{N \to \infty} \frac{-\log B_{h,\beta}(N)}{N}$$

*of $B_{h,\beta}$ exists in $[0, \varphi_\beta(1)]$ and is continuous as function on $\mathbb{R}^+ \times \mathbb{R}^+$, and*

(3.5) $$B_{h,\beta}(N) \leq e^{-m_B(h,\beta)N}$$

*is valid for all $N \in \mathbb{N}$. For $h > 0$ and $\beta < \beta_c(h)$, we further have*

(3.6) $$m_B(h,\beta) = m_G(h,\beta).$$

*Moreover, for $h > 0$ and $\beta_0 < \beta_c(h)$, we also have*

(3.7) $$B_{h,\beta}(N) \geq \frac{1}{K_{h,\beta_0}} e^{-m_B(h,\beta)N}$$

*for all $N \in \mathbb{N}$ and $\beta \leq \beta_0$, where $K_{h,\beta_0} < \infty$ is defined in (3.4).*

PROOF. By (3.2) in Lemma 3.4, the existence of $m_B$ in $\mathbb{R} \cup \{-\infty\}$ and the estimate in (3.5) are consequences of the subadditive limit lemma applied to the sequence $(-\log B_{h,\beta}(N))_{N \in \mathbb{N}}$. The lower bound for $m_B$ is obvious and continuity is obtained by similar arguments as in the proof of Proposition 2.7. The upper bound $\varphi_\beta(1)$ for $m_B$ follows from the lower estimate in Remark 3.3 once we have shown (3.6).

Now, suppose that $h > 0$ and $\beta < \beta_c(h)$. For the proof of (3.6), recall that Theorem A states that

(3.8) $$m_G(h,\beta) = \lambda_h - \lambda_{\bar{h}},$$

where $\lambda_{h'} = \log E_0[\exp(h' \cdot S_1(1))]$ for $h' \geq 0$ and where $\bar{h} = \bar{h}(h,\beta) > 0$ is determined by

(3.9) $$\overline{m}_G(\bar{h},\beta) = h - \bar{h},$$

with $\overline{m}_G$ having been introduced in Section 2.1. By the definition of $P_h$ in (1.1), we consequently have

(3.10)
$$b_{h,\beta}(L;N)e^{m_G(h,\beta)N} = b_{0,\beta}(L;N)e^{-\lambda_{\bar{h}}N+hL}$$
$$= b_{\bar{h},\beta}(L;N)e^{\overline{m}_G(\bar{h},\beta)L}$$

and therefore, since $\overline{m}_G(\bar{h},\beta) = \overline{m}_B(\bar{h},\beta)$ by Corollary 2.9,

$$\sum_{N=1}^{\infty} B_{h,\beta}(N)e^{m_G(h,\beta)N} = \sum_{L=1}^{\infty} \overline{B}_{\bar{h},\beta}(L)e^{\overline{m}_B(\bar{h},\beta)L} = \infty,$$

where the second equality follows from the lower estimate for $\overline{B}_{\bar{h},\beta}(L)$ in Proposition 2.7. As a consequence, the mass $m_B(h,\beta)$ cannot be greater than $m_G(h,\beta)$. Since the inverted estimate is obvious, this proves (3.6).

It remains to show (3.7). By the definition of $P_h$, it is plain that

$$P_h[S_1(N) = 0] \leq e^{-\lambda_h N}.$$

In the sub-ballistic regime, by (3.6) and (3.8), we further have that

$$m_B(h,\beta) = m_G(h,\beta) < \lambda_h.$$

Therefore, and by obvious monotonicity, we get

$$K_{h,\beta} \leq K_{h,\beta_0} < \infty$$

for all $\beta \leq \beta_0$ so that (3.7) now follows from (3.3) and the subadditive limit lemma applied to $\log(K_{h,\beta}B_{h,\beta}(N))_{N\in\mathbb{N}}$. □

COROLLARY 3.6. *For any $h > 0$ and $\beta_0 < \beta_c(h)$, we have*

$$G_{h,\beta}(N) \geq \frac{1}{K_{h,\beta_0}} e^{-m_G(h,\beta)N}$$

*for all $N \in \mathbb{N}$ and $\beta \leq \beta_0$, where $K_{h,\beta_0} < \infty$ is defined in (3.4).*

PROOF. The corollary follows from (3.6) and (3.7). □

3.2. *Exponential gap and analyticity.* As in Section 2.2, for fixed $p \in \mathbb{N}$, we consider independent copies

$$\mathcal{S}^j = (S^j(n))_{n\in\mathbb{N}_0}, \qquad j = 1,\ldots,p,$$

of the random walk $\mathcal{S}$, defined on the probability space $(\Omega^p, \mathcal{F}^{\otimes p}, P_h)$. The random process $\mathcal{S}^{(p)} = (S^{(p)}(n))_{n\in\mathbb{N}_0^p}$ with values in $(\mathbb{Z}^d)^p$ is given by

$$S^{(p)}(n) = (S^1(n^1),\ldots,S^p(n^p)), \qquad n = (n^1,\ldots,n^p) \in \mathbb{N}_0^p.$$

The potential $\Phi_\beta^{(p)}$ for $\mathcal{S}^{(p)}$ was introduced as

$$\Phi_\beta^{(p)}(N) = \sum_{j=1}^p \Phi_\beta^j(N^j), \qquad N = (N^1, \ldots, N^p) \in \mathbb{N}_0^p,$$

where $\Phi_\beta^1, \ldots, \Phi_\beta^p$ are the corresponding potentials associated with the single random walks $\mathcal{S}^1, \ldots, \mathcal{S}^p$.

For $\omega \in \Omega^p$ and $M, N \in \mathbb{N}_0^p$ with $M \leq N$ componentwise, recall that the path

$$\mathcal{S}^{(p)}[M, N](\omega) = (S^{(p)}(n, \omega))_{n \in \mathbb{N}_0 : M \leq n \leq N}$$

is called a bridge if and only if the paths

$$\mathcal{S}^j[M^j, N^j](\omega), \qquad j = 1, \ldots, p,$$

are bridges as single random walks, starting in the hyperplane $\mathcal{H}_{S_1^1(M^1, \omega)}$ and ending in the hyperplane $\mathcal{H}_{S_1^1(N^1, \omega)}$. Moreover, for $M < N$, a bridge $\mathcal{S}^{(p)}[M, N](\omega)$ is called irreducible if and only if $S_1^1(N^1, \omega)$ is its only (common) breaking point.

For $h, \beta \geq 0$, $m \in \mathbb{N}_0^p$ and $n \in \mathbb{N}^p$ with $m < n$, we now define

$$\mathrm{Ir}^p(n) \stackrel{\mathrm{def}}{=} \{\mathcal{S}^{(p)}[n] \text{ is an irreducible bridge}\},$$

$$\mathrm{Ir}^p(m, n) \stackrel{\mathrm{def}}{=} \{\mathcal{S}^{(p)}[m, n] \text{ is an irreducible bridge}\}$$

and

$$\Lambda_{h,\beta}^p(n) \stackrel{\mathrm{def}}{=} E_h[\exp(-\Phi_\beta^{(p)}(n)); \mathrm{Ir}^p(n)].$$

Observe that we have

$$\mathrm{Ir}^p(n) = \bigcup_{L \in \mathbb{N}} \mathrm{ir}^p(L; n) \quad \text{and} \quad \Lambda_{h,\beta}^p(n) = \sum_{L \in \mathbb{N}} \lambda_{h,\beta}^p(L; n),$$

where the union is of disjoint sets and where

$$\mathrm{ir}^p(L; n) = \{\mathcal{S}^{(p)}[n] \text{ is an irreducible bridge of span } L\},$$

$$\lambda_{h,\beta}^p(L; n) = E_h[\exp(\Phi_\beta^{(p)}(n)); \mathrm{ir}^p(L; n)]$$

were introduced in Section 2.2.

By means of Theorem A, we are able to transfer the mass gap for irreducible bridges from Section 2.3 to the present $N$-step setting, detecting an exponential gap between the long-time behavior of irreducible bridges and arbitrary bridges.

THEOREM 3.7. *Suppose $h > 0$ and $\beta_0 < \beta_c(h)$. There is then some $\gamma > 0$ such that*

$$\Lambda^p_{h,\beta}(n) e^{m_B(h,\beta) \sum_{j=1}^p n^j} \leq \frac{1}{\gamma} e^{-\gamma \sum_{j=1}^p n^j} \tag{3.11}$$

*for all $n = (n^1, \ldots, n^p) \in \mathbb{N}^p$ and $\beta \leq \beta_0$. Moreover, we have*

$$\sum_{n \in \mathbb{N}^p} \Lambda^p_{h,\beta}(n) e^{m_B(h,\beta) \sum_{j=1}^p n^j} = 1 \tag{3.12}$$

*for all $\beta < \beta_c(h)$.*

PROOF. We first deal with (3.12). From (3.8) and (3.9), and by letting $N \stackrel{\text{def}}{=} \sum_{j=1}^p n^j$, we obtain

$$\begin{aligned}
\lambda^p_{h,\beta}(L;n) e^{m_G(h,\beta)N} &= \lambda^p_{0,\beta}(L;n) e^{-\lambda_{\bar h} N + phL} \\
&= \lambda^p_{\bar h,\beta}(L;n) e^{p \overline{m}_G(\bar h,\beta)L}
\end{aligned} \tag{3.13}$$

by an analogous argument as for (3.10). In the ballistic regime, since $\overline{m}_G(\bar h, \beta) = \overline{m}_B(\bar h, \beta)$ and $p\overline{m}_B(\bar h, \beta) = \overline{m}^p_B(\bar h, \beta)$, we have

$$\sum_{n \in \mathbb{N}^p} \Lambda^p_{h,\beta}(n) e^{m_G(h,\beta)N} = \sum_{L \in \mathbb{N}} \overline{\Lambda}^p_{\bar h,\beta}(L) e^{\overline{m}^p_B(\bar h,\beta)L} = 1,$$

where the second equality is part of Lemma 2.15.

In order to achieve an exponential gap, observe that for any $\delta > 0$, again by (3.13), we have

$$\begin{aligned}
\Lambda^p_{h,\beta}(n) e^{m_G(h,\beta)N} &= \sum_{L \leq \delta N} \lambda^p_{h,\beta}(L;n) e^{m_G(h,\beta)N} + \sum_{L > \delta N} \lambda^p_{\bar h,\beta}(L,n) e^{\overline{m}^p_B(\bar h,\beta)L} \\
&\leq \sum_{L \leq \delta N} \lambda^p_{h,0}(L;n) e^{m_G(h,\beta)N} + \sum_{L > \delta N} \overline{\Lambda}^p_{\bar h,\beta}(L) e^{\overline{m}^p_B(\bar h,\beta)L}.
\end{aligned}$$

From the independence of $\mathcal{S}^1, \ldots, \mathcal{S}^2$, the definition of $P_h$ and (3.8) again, we obtain

$$\sum_{L \leq \delta N} \lambda^p_{h,0}(L;n) e^{m_G(h,\beta)N} \leq \prod_{j=1}^p P_h[S^j_1(n^j) \leq \delta N] e^{m_G(h,\beta)n^j} \leq e^{(h\delta p - \lambda_{\bar h})N}.$$

Moreover, by Lemma 2.17, we have

$$\begin{aligned}
\sum_{L > \delta N} \overline{\Lambda}^p_{\bar h,\beta}(L) e^{\overline{m}^p_B(\bar h,\beta)L} &\leq \frac{1}{p} e^{2(\varphi_\beta(1) + \lambda_{\bar h})} \sum_{L > \delta N} e^{-(\overline{m}^p_\Lambda(\bar h,\beta) - \overline{m}^p_B(\bar h,\beta))L} \\
&\leq \frac{1}{p} e^{2(\varphi_\beta(1) + \lambda_{\bar h})} \frac{e^{-\delta(\overline{m}^p_\Lambda(\bar h,\beta) - \overline{m}^p_B(\bar h,\beta))N}}{1 - e^{-(\overline{m}^p_\Lambda(\bar h,\beta) - \overline{m}^p_B(\bar h,\beta))}}.
\end{aligned}$$

Since $\lambda_{\bar{h}} > 0$ by Theorem A and $\overline{m}_\Lambda^p(\bar{h}, \beta) > \overline{m}_B^p(\bar{h}, \beta)$ by Theorem 2.18, this proves (3.11) for a single $\beta \leq \beta_c(h)$.

In order to find an uniform estimate, observe that

$$\overline{m}_G(h', \beta) + h' = \lim_{L \to \infty} \frac{-\log(\sum_{N \in \mathbb{N}} E_0[e^{-\Phi_\beta(N) - \lambda_{h'} N}; \{S_1(N) = L\}])}{L}$$

is increasing in both variables. Since $\bar{h} = \bar{h}(h, \beta)$ fulfills $\overline{m}_B(\bar{h}, \beta) + \bar{h} = h$, we thus have $\lambda_{\bar{h}(h,\beta)} \leq \lambda_{\bar{h}(h,\beta_0)}$ for $\beta \leq \beta_0$. The existence of a uniform bound now follows from the continuity of $\overline{m}_B^p$ and $\overline{m}_\Lambda^p$ on $\mathbb{R}^+ \times \mathbb{R}^+$, which we established in Section 2. □

By known arguments using the analytic implicit function theorem (see page 329 of [9]), we obtain the following first consequence of Theorem 3.7.

COROLLARY 3.8. *Suppose that the function $\varphi$, introduced at the beginning of Section 2, is analytic on $(0, \infty)$. The mass $m_G$ is then analytic on the open set $\{(h, \beta) \in (0, \infty) : \beta < \beta_c(h)\}$.*

REMARK. By dominated convergence, it is obvious that the function

$$\varphi(t) = -\log \mathbb{E} \exp(-tV_x)$$

is analytic on $(0, \infty)$. Therefore, by Proposition 3.2, Corollary 3.6 and now Corollary 3.8, we have completed the proof of Theorem B in Section 1.

3.3. *Restricted path intersections.* In this section, we investigate finite random walks $\mathcal{S}[N]$ with restricting assumptions on

$$\|\ell(N)\|_2 \stackrel{\text{def}}{=} \left( \sum_{z \in \mathbb{Z}^d} \ell_z(N)^2 \right)^{1/2}, \qquad N \in \mathbb{N}.$$

More precisely, for $h, \beta \geq 0$, $N \in \mathbb{N}$ and $k \in \mathbb{R}^+$, we define

$$G_{h,\beta}^{\leq k}(N) \stackrel{\text{def}}{=} E_h[\exp(-\Phi_\beta(N)); \{\|\ell(N)\|_2^2 \leq kN\}]$$

and we want to show that in the ballistic regime, such restrictions have no crucial effect on the long-time behavior.

To this end, we first investigate restricted bridges. For $h, \beta \geq 0$, $N \in \mathbb{N}$ and $k \in \mathbb{R}^+$, we define

$$B_{h,\beta}^{>k}(N) \stackrel{\text{def}}{=} E_h[\exp(-\Phi_\beta(N)); \{\|\ell(N)\|_2^2 > kN\} \cap \text{Br}(N)].$$

PROPOSITION 3.9. *For any $h > 0$ and $\beta_0 < \beta_c(h)$, we have*

$$\lim_{k \to \infty} \sup_{\beta \leq \beta_0} \sup_{N \in \mathbb{N}} B_{h,\beta}^{>k}(N) e^{m_G(h,\beta)N} = 0.$$

PROOF. Let $(\tau_i, \xi_i)_{i \in \mathbb{N}}$, be a sequence of independent, identically distributed random vectors with distribution

$$P_{h,\beta}[(\tau_i, \xi_i) = (n, x)] \stackrel{\text{def}}{=} e^{m_G(h,\beta)n} E_h[e^{-\Phi_\beta(n)}; \{\|\ell(n)\|_2^2 = x\} \cap \text{Ir}(n)]$$

for $n, x \in \mathbb{N}$ and $\beta < \beta_c(h)$. In fact, $P_{h,\beta}$ is a probability distribution by (3.12). The expectation with respect to $P_{h,\beta}$ will be denoted by $E_{h,\beta}$.

Now, assume that $\omega \in \text{Br}(N)$ and let $m \in \mathbb{N}$ be the number of breaking points for $\mathcal{S}[N](\omega)$. There then exist unique times $0 = n_0 < \cdots < n_m = N$ such that $\mathcal{S}[n_0, n_1](\omega), \ldots, \mathcal{S}[n_{m-1}, n_m](\omega)$ are irreducible bridges. We thus have

$$\text{Br}(N) = \bigcup_{m=1}^{N} \bigcup_{0 = n_0 < \cdots < n_m = N} \bigcap_{i=1}^{m} \text{Ir}(n_{i-1}, n_i),$$

where the union is of disjoint sets. For every $\omega \in \bigcap_{i=1}^{m} \text{Ir}(n_{i-1}, n_i)$, we further have

$$\|\ell(N, \omega)\|_2^2 = \sum_{i=1}^{m} \sum_{z \in \mathbb{Z}^d} \ell_z(n_{i-1}, n_i, \omega)^2,$$

$$\Phi_\beta(N, \omega) = \sum_{i=1}^{m} \Phi_\beta(n_{i-1}, n_i, \omega),$$

where the second equation goes back to Lemma 2.1(b). As a consequence, by the Markov property of $\mathcal{S}$, we obtain that $B_{h,\beta}^{>k}(N) e^{m_G(h,\beta)N}$ equals

$$\sum_{m=1}^{N} \sum_{x_1 + \cdots + x_m > kN} \sum_{0 = n_0 < \cdots < n_m = N}$$
$$\prod_{i=1}^{m} e^{m_G(h,\beta)(n_i - n_{i-1})} E_h[e^{-\Phi_\beta(n_i - n_{i-1})} 1_{\text{Ir}(n_i - n_{i-1})} 1_{\{\|\ell(n_i - n_{i-1})\|_2^2 = x_i\}}]$$

$$= \sum_{m=1}^{N} P_{h,\beta}\left[\sum_{1 \leq i \leq m} \tau_i = N, \sum_{1 \leq i \leq m} \xi_i > Nk\right]$$

$$\leq N P_{h,\beta}\left[\sum_{1 \leq i \leq N} \xi_i > Nk\right]$$

for all $k \in \mathbb{R}^+$.

Next, observe that for any $n \in \mathbb{N}$, we have

$$\|\ell(n)\|_2 \leq \sum_{z \in \mathbb{Z}^d} \ell_z(n) = n$$

such that, by the exponential gap in Theorem 3.7, the moments

$$E_{h,\beta}[\xi_1^m] \leq \sum_{n \in \mathbb{N}} n^{2m} \Lambda_{h,\beta}(n) e^{m_G(h,\beta)n}, \qquad m \in \mathbb{N},$$

are finite and continuous in $\beta < \beta_c(h)$. For any $k > E_{h,\beta}[\xi_1]$, by the Chebyshev inequality and the independence of $(\xi_i)_{i \in \mathbb{N}}$, we moreover have

$$B_{\beta,h}^{>k}(N) e^{m_G(h,\beta)N} \leq N P_{h,\beta}\left[\sum_{1 \leq i \leq N} |\xi_i - E_{h,\beta}[\xi_1]| > N|k - E_{h,\beta}[\xi_1]|\right]$$

$$\leq \frac{1}{N(k - E_{h,\beta}[\xi_1])^2} E_{h,\beta}\left[\left(\sum_{1 \leq i \leq N} (\xi_i - E_{h,\beta}[\xi_1])\right)^2\right]$$

$$= \frac{1}{(k - E_{h,\beta}[\xi_1])^2} E_{h,\beta}[(\xi_1 - E_{h,\beta}[\xi_1])^2].$$

The proposition now follows from the continuity of the first and second moments of $\xi_1$ in $\beta < \beta_c(h)$. $\square$

By means of Proposition 3.9, we are now able to prove Lemma E from Section 1, formulated in the present setting of a generalized potential.

COROLLARY 3.10. *For any $h > 0$, $\beta_0 < \beta_c(h)$ and $\varepsilon < 1/K_{h,\beta_0}$, there exists $k_\varepsilon < \infty$ such that*

(3.14) $$G_{h,\beta}^{\leq k_\varepsilon}(N) \geq \varepsilon G_{h,\beta}(N)$$

*for all $N \in \mathbb{N}$ and $\beta \leq \beta_0$.*

PROOF. For any $k \in \mathbb{R}^+$ and $N \in \mathbb{N}$, we obviously have

$$G_{h,\beta}^{\leq k}(N) \geq B_{h,\beta}(N) - B_{h,\beta}^{>k}(N).$$

Moreover, by Proposition 3.2 and Corollary 3.6, we know that

$$G_{h,\beta}(N) \leq e^{-m_G(h,\beta)N} \leq K_{h,\beta_0} B_{h,\beta}(N)$$

for all $N \in \mathbb{N}$ and $\beta \leq \beta_0$. The corollary thus follows from Lemma 3.9. $\square$

**4. Coupled path potential.** We consider two independent copies $\mathcal{S}^1$ and $\mathcal{S}^2$ of the random walk $\mathcal{S}$ with drift $h \geq 0$ and starting condition

$$P_{h,y^1,y^2}[S^1(0) = (0, y^1), S^2(0) = (0, y^2)] = 1$$

for $y^1, y^2 \in \mathbb{Z}^{d-1}$. Expectations with respect to $P_{h,y^1,y^2}$ are denoted by $E_{h,y^1,y^2}$, where the indices $y^1, y^2$ are left off when $\mathcal{S}^1$ and $\mathcal{S}^2$ start at the origin.

As in Section 2.2, we compose from the two random walks the random process

$$\mathcal{S}^{(2)} = (S^{(2)}(n))_{n \in \mathbb{N}_0^2}$$

on $\mathbb{Z}^d \times \mathbb{Z}^d$, where $S^{(2)}(n)$, for $n = (n^1, n^2)$, is given by

$$S^{(2)}(n) = (S^1(n^1), S^2(n^2)).$$

For any parameter $\beta \geq 0$, motivated by the heuristic picture of Theorem D in Section 1, we introduce a coupled path potential $\widetilde{\Phi}_\beta^{(2)}$ for the process $\mathcal{S}^{(2)}$. For $M = (M^1, M^2)$ and $N = (N^1, N^2) \in \mathbb{N}_0$ with $M \leq N$ componentwise, we define

$$\widetilde{\Phi}_\beta^{(2)}(N) \stackrel{\text{def}}{=} \sum_{x \in \mathbb{Z}^d} \varphi_\beta(\ell_x^1(N^1) + \ell_x^2(N^2)),$$

$$\widetilde{\Phi}_\beta^{(2)}(M, N) \stackrel{\text{def}}{=} \sum_{x \in \mathbb{Z}^d} \varphi_\beta(\ell_x^1(M^1, N^1) + \ell_x^2(M^2, N^2)),$$

where $\varphi_\beta$ was introduced at the beginning of Section 2 and where

$$\ell_x^j(N^j) = \sum_{n=1}^{N^j} 1_{\{S^j(n)=x\}}, \qquad \ell_x^j(M, N) = \sum_{n=M^j+1}^{N^j} 1_{\{S^j(n)=x\}}$$

denote the number of visits to the site $x \in \mathbb{Z}^d$ by the random walk $\mathcal{S}^j[1, N]$, respectively $\mathcal{S}^j[M+1, N]$. Observe that in contrast to the path potential $\Phi_\beta^{(2)}$, the distribution of this coupled potential $\widetilde{\Phi}_\beta^{(2)}$ depends on the starting sites $y^1, y^2 \in \mathbb{Z}^{d-1}$ of the random walks $\mathcal{S}^1$ and $\mathcal{S}^2$ (i.e., on $\|y^2 - y^1\|$).

For $h, \beta \geq 0$, $y^1, y^2 \in \mathbb{Z}^{d-1}$ and $N \in \mathbb{N}^2$, we define

$$\widetilde{G}_{h,\beta,y^1,y^2}^2(N) \stackrel{\text{def}}{=} E_{h,y^1,y^2}[\exp(-\widetilde{\Phi}_\beta^{(2)}(N))].$$

The aim of this section is to prove the following "second moment"-type estimate on $\widetilde{G}_{h,\beta,y^1,y^2}^2$.

THEOREM 4.1. *Suppose that $d \geq 4$ and $h > 0$. There then exist $\beta_0 > 0$ and $K_{\text{s.m.}} < \infty$ such that*

$$\widetilde{G}_{h,\beta,y^1,y^2}^2(N) \leq K_{\text{s.m.}} G_{h,\beta}(N^1) G_{h,\beta}(N^2)$$

*for all $N = (N^1, N^2) \in \mathbb{N}^2$, $y^1, y^2 \in \mathbb{Z}^{d-1}$ and $\beta \leq \beta_0$.*

COINCIDENCE OF LYAPUNOV EXPONENTS 49

REMARK. In the particular case $\varphi(t) = -\log \mathbb{E} e^{-tV_x}$, we have

$$\widetilde{G}^2_{h,\beta,y^1,y^2}(N) = \mathbb{E} Z^{\mathrm{qu}}_{\mathbb{V},h,\beta,N^1,y^1} Z^{\mathrm{qu}}_{\mathbb{V},h,\beta,N^2,y^2},$$

where $Z^{\mathrm{qu}}_{\mathbb{V},h,\beta,N^j,y^j}$ denotes the quenched partition function from Section 1, but here with starting condition $P_h[S^j(0) = (0, y^j)] = 1$ for $j = 1, 2$. Theorem D is thus a special case of Theorem 4.1.

In order to establish Theorem 4.1, we investigate bridges under $\widetilde{\Phi}^{(2)}_\beta$. For $h, \beta \geq 0$, $y^1, y^2 \in \mathbb{Z}^{d-1}$ and $M, N \in \mathbb{N}_0^2$ with $M \leq N$, we set

$$\mathrm{Br}^2(N) \stackrel{\mathrm{def}}{=} \{\mathcal{S}^{(2)}[N] \text{ is a bridge}\},$$

$$\mathrm{Br}^2(M,N) \stackrel{\mathrm{def}}{=} \{\mathcal{S}^{(2)}[M,N] \text{ is a bridge}\}$$

and

$$\widetilde{B}^2_{h,\beta,y^1,y^2}(N) \stackrel{\mathrm{def}}{=} E_{h,y^1,y^2}[\exp(-\widetilde{\Phi}^{(2)}_\beta(N)); \mathrm{Br}^2(N)].$$

We want to divide the bridges into irreducible "strips" which may then be treated by renewal techniques. To this end, for $h > 0$ and $\beta < \beta_c(h)$, let $(\tau_i^1, \tau_i^2, \eta_i^1, \eta_i^2, \zeta_i)_{i \in \mathbb{N}_0}$ be a Markov chain with transition probabilities

$$P_{h,\beta}[(\tau_{i+1}^1, \tau_{i+1}^2, \eta_{i+1}^1, \eta_{i+1}^2, \zeta_{i+1})$$
$$= (n^1, n^2, y^1, y^2, z) | (\tau_i^1, \tau_i^2, \eta_i^1, \eta_i^2, \zeta_i)]$$
$$\stackrel{\mathrm{def}}{=} e^{m_G(h,\beta)(n^1+n^2)} E_{h,\eta_i^1,\eta_i^2}[e^{-\Phi^{(2)}_\beta(n)} 1_{\mathrm{Ir}^2(n)} 1_{\{S^1(n^1)=(S_1^1(n^1),y^1)\}}$$
$$\times 1_{\{S^2(n^2)=(S_1^2(n^2),y^2)\}} 1_{\{L(n)=z\}}]$$

for $n = (n^1, n^2) \in \mathbb{N}^2$, $y^1, y^2 \in \mathbb{Z}^{d-1}$ and $z \in \mathbb{N}_0$, and with

$$L(n) \stackrel{\mathrm{def}}{=} \sum_{x \in R^1(n^1) \cap R^2(n^2)} \ell_x^1(n^1) + \ell_x^2(n^2),$$

where $R^j(n^j) \stackrel{\mathrm{def}}{=} \{x \in \mathbb{Z}^d : \ell_x^j(n^j) > 0\}$ for $j = 1, 2$. In fact, $P_{h,\beta}$ is a probability distribution by (3.12). We will write $P_{h,\beta,y^1,y^2}$ to indicate the starting condition

$$P_{h,\beta,y^1,y^2}[\eta_0^1 = y^1, \eta_0^2 = y^2] = 1,$$

and expectations with respect to $P_{h,\beta,y^1,y^2}$ are denoted by $E_{h,\beta,y^1,y^2}$. Again, if the start is at the origin, the indices are omitted.

In order to bound $\widetilde{\Phi}_\beta$ within an irreducible strip, observe that

(4.1) $\qquad \widetilde{\Psi}^{(2)}_\beta(n) \stackrel{\mathrm{def}}{=} \Phi^{(2)}_\beta(n) - \widetilde{\Phi}^{(2)}_\beta(n) \leq \varphi_\beta(1) L(n)$

for all $n = (n^1, n^2) \in \mathbb{N}^2$. It is thus convenient to define $\sigma_0 \stackrel{\text{def}}{=} 0$ and

$$\sigma_k \stackrel{\text{def}}{=} \min\{i > \sigma_{k-1} : \zeta_i > 2 \cdot 1_{\{\eta_{i-1}^1 = \eta_{i-1}^2\}} 1_{\{\eta_i^1 \neq \eta_i^2\}}\}, \qquad k \in \mathbb{N},$$

as well as $\rho_0 \stackrel{\text{def}}{=} 1_{\{\eta_0^1 = \eta_0^2\}}$ and

$$\rho_i \stackrel{\text{def}}{=} \begin{cases} \zeta_i + 2 \cdot 1_{\{\eta_i^1 = \eta_i^2\}}, & \text{if } \exists k \in \mathbb{N} \text{ with } \sigma_k = i,\ i \in \mathbb{N}, \\ 0, & \text{otherwise.} \end{cases}$$

That means that we want to know in which strips we have path intersections of $\mathcal{S}^1$ and $\mathcal{S}^2$, $\sigma_k$ denoting the $k$th of these strips. However, if $\mathcal{S}^1$ and $\mathcal{S}^2$ enter a strip from a common site, then they must go "forward" at the first step and they consequently intersect for a first time. If this strip does not contain any further intersections (and does not consist of only one step), then it is not considered in the definition of $\sigma_k$. For such a strip, the contribution to $\widetilde{\Psi}_\beta$ is anticipated in the previous strip by the summand $1_{\{\eta_i^1 = \eta_i^2\}}$ in the definition of $\rho_i$. This special treatment of such strips is necessary to have $P_{h,\beta}[\sigma_1 > 1] > 0$.

For $m \in \mathbb{N}_0$, let $T(m) = (T^1(m), T^2(m))$ now be given by

$$T^j(m) \stackrel{\text{def}}{=} \sum_{1 \leq i \leq m} \tau_i^j$$

for $j = 1, 2$. The conclusion of the above comments on the definitions of $\sigma_k$ and $\rho_i$ is the following.

LEMMA 4.2. *For any $h > 0$ and $\beta < \beta_c(h)$, we have*

$$\widetilde{B}_{h,\beta,y_0^1,y_0^2}^2(N) e^{m_G(h,\beta)(N^1+N^2)} \leq \sum_{m \in \mathbb{N}} E_{h,\beta,y_0^1,y_0^2}[e^{\varphi_\beta(1) \sum_{i=0}^m \rho_i} 1_{\{T(m) = N\}}]$$

*for all $N = (N^1, N^2) \in \mathbb{N}^2$ and $y_0^1, y_0^2 \in \mathbb{Z}^{d-1}$.*

PROOF. By similar arguments as in the proof of Lemma 3.9, we obtain

$$\text{Br}^2(N) = \bigcup_{m \in \mathbb{N}} \bigcup_{\substack{n_0,\ldots,n_m \in \mathbb{N}_0^2 : \\ 0 = n_0 < \cdots < n_m = N}} \bigcap_{i=1}^m \text{Ir}^2(n_{i-1}, n_i),$$

where the union is of disjoint sets and where $m \in \mathbb{N}$ represents the number of breaking points for the corresponding $N$-step bridge. For $\omega \in \bigcap_{i=1}^m \text{Ir}^2(n_{i-1}, n_i)$, we further have

$$\widetilde{\Phi}_\beta^{(2)}(n, \omega) = \sum_{i=1}^m \Phi_\beta^{(2)}(n_{i-1}, n_i, \omega) - \widetilde{\Psi}_\beta^{(2)}(n_{i-1}, n_i, \omega).$$

Therefore, by renewing the random walk $\mathcal{S}^{(2)}$ at times $n_1, \ldots, n_{m-1}$, we obtain that $\widetilde{B}^2_{h,\beta,y_0^1,y_0^2}(N) e^{m_G(h,\beta)(N^1+N^2)}$ is equal to

$$\sum_{m=1}^{\infty} \sum_{\substack{n_0,\ldots,n_m \in \mathbb{N}_0^2: \\ 0=n_0<\cdots<n_m=N}} \sum_{y_1^1,y_1^2,\ldots,y_m^1,y_m^2 \in \mathbb{Z}^{d-1}} \sum_{z_1,\ldots,z_m \in \mathbb{N}_0} \prod_{i=1}^m e^{m_G(h,\beta)(n_i^1-n_{i-1}^1+n_i^2-n_{i-1}^2)}$$

$$\times E_{h,y_{i-1}^1,y_{i-1}^2}[e^{-\Phi_\beta^{(2)}(n_i-n_{i-1})+\widetilde{\Psi}_\beta^{(2)}(n_i-n_{i-1})} 1_{\mathrm{Ir}^2(n_i-n_{i-1})}$$

$$\times 1_{\{S^j(n_i^j-n_{i-1}^j)=(S_j^1(n_i^j-n_{i-1}^j),y_i^j) \text{ for } j=1,2\}} 1_{\{L(n_i-n_{i-1})=z_i\}}]$$

$$\leq \sum_{m=1}^{\infty} E_{h,\beta,y_0^1,y_0^2}[e^{\varphi_\beta(1)\sum_{i=0}^m \rho_i} 1_{\{T^j=n^j \text{ for } j=1,2\}}],$$

where we assume $n_i = (n_i^1, n_i^2)$ for $i \in \mathbb{N}$ and where (4.1) is used. $\square$

The value of this renewal formalism is substantiated by the following estimate.

PROPOSITION 4.3. *For any $h > 0$ and $\beta_0 < \beta_c(h)$, and with $K_{h,\beta_0} < \infty$ being defined in (3.4), we have*

$$\frac{\widetilde{G}^2_{h,\beta,y^1,y^2}(N)}{G_{h,\beta}(N^1) G_{h,\beta}(N^2)} \leq 4 K_{h,\beta_0}^4 \sup_{y_0^1,y_0^2} \sum_{k=0}^{\infty} E_{h,\beta,y_0^1,y_0^2}[e^{\varphi_\beta(1)\sum_{i=0}^k \rho_{\sigma_i}} 1_{\{\sigma_k<\infty\}}]$$

*for all $N = (N^1, N^2) \in \mathbb{N}^2$, $y^1, y^2 \in \mathbb{Z}^{d-1}$ and $\beta \leq \beta_0$.*

PROOF. For $j \in \{1,2\}$, by the definition of $\widetilde{\Phi}_\beta$ and the monotonicity of $\varphi_\beta$, we obviously have

$$\widetilde{\Phi}_\beta^{(2)}(N) \geq \Phi_\beta^j(N^j),$$

where $\Phi_\beta^j$ is the single path potential associated with the random walk $\mathcal{S}^j$. From the independence of $\mathcal{S}^1$ and $\mathcal{S}^2$ thus follows

$$E_{h,y^1,y^2}[e^{-\widetilde{\Phi}_\beta^{(2)}(N)}; \{S_1^1(N^1) \leq 0 \text{ or } S_1^2(N^2) \leq 0\}]$$

$$\leq E_h[e^{-\Phi_\beta^1(N^1)}; \{S_1^2(N^2) \leq 0\}] + E_h[e^{-\Phi_\beta^2(N^2)}; \{S_1^1(N^1) \leq 0\}]$$

$$= G_{h,\beta}(N^1) P_h[S^2(N^2) \leq 0] + G_{h,\beta}(N^2) P_h[S^1(N^1) \leq 0]$$

$$\leq 2 K_{h,\beta_0} G_{h,\beta}(N^1) G_{h,\beta}(N^2)$$

for all $N^1, N^2 \in \mathbb{N}$ and $\beta \leq \beta_0$, where the last step goes back to (3.8) and Corollary 3.6.

We now consider the more complicated case of positive first components. Since $\mathcal{S}^1$ and $\mathcal{S}^2$ are exchangeable, it suffices to consider

$$\omega \in \{0 < S_1^1(N^1) \leq S_1^2(N^2)\}.$$

In that case, $\omega$ is also an element of the union

$$\bigcup_{0<\bar{N}^2\leq N^2} \bigcup_{0\leq m_1^1<m_2^1\leq N^2} \bigcup_{0\leq m_1^2<m_2^2\leq \bar{N}^2}$$

$$\{S_1^1(m_1^1) = S_1^2(m_1^2) = 0\} \cap \operatorname{Br}^2(m_1, m_2) \cap \{S_1^1(N^1) = S_1^1(m_2^1)\}$$
$$\cap \{S_1^2(\bar{N}^2) = S_1^2(m_2^2)\} \cap \{S_1^2(\bar{N}^2) < S_1^2(\nu) \text{ for } \bar{N}^2 < \nu \leq N^2\},$$

where $m_1 \stackrel{\text{def}}{=} (m_1^1, m_1^2)$ and $m_2 \stackrel{\text{def}}{=} (m_2^1, m_2^2)$ and where $m_1^j$ and $m_2^j$ may be chosen as the times of the last visit of $\mathcal{S}^j[N^j](\omega)$ to the hyperplane $\mathcal{H}_0$, respectively the first visit of $\mathcal{S}^j[m_1^j, N^j](\omega)$ to the hyperplane $\mathcal{H}_{S_1^1(N^1,\omega)}$, and where $\bar{N}^2$ is the last-visit time of $\mathcal{S}^2[N^2](\omega)$ to $\mathcal{H}_{S_1^1(N^1,\omega)}$. Moreover, we have

$$\widetilde{\Phi}_\beta^{(2)}(N, \omega) \geq \widetilde{\Phi}_\beta^{(2)}(m_1, m_2, \omega) + \Phi_\beta^2(\bar{N}^2, N^2, \omega)$$

for the corresponding $m_1$, $m_2$ and $\bar{N}^2$. Therefore, with the notation

$$G_{h,\beta}^+(k) \stackrel{\text{def}}{=} E_h[e^{-\Phi_\beta^2(N)}; \{0 < S_1^2(\nu) \text{ for } \nu = 1, \ldots, k\}], \qquad k \in \mathbb{N}_0,$$

and by renewing $\mathcal{S}^{(2)}$ at times $m_1, m_2$ and $\bar{N}^2$, we obtain

$$E_{h,y^1,y^2}[e^{-\widetilde{\Phi}_\beta^{(2)}(N)} 1_{\{0<S_1^1(N^1)\leq S_1^2(N^2)\}}]$$

$$\leq \sum_{0<\bar{N}^2\leq N^2} \sum_{0\leq m_1^1<m_2^1\leq N^1} \sum_{0\leq m_1^2<m_2^2\leq \bar{N}^2} \sum_{y_0^1,y_0^2\in\mathbb{Z}^{d-1}}$$

$$E_{h,y^1,y^2}\big[1_{\{S^1(m_1^1)=(0,y_0^1)\}} 1_{\{S^2(m_1^2)=(0,y_0^2)\}} e^{-\widetilde{\Phi}_\beta^{(2)}(m_1,m_2)}$$
$$\times 1_{\operatorname{Br}^2(m_1,m_2)} 1_{\{S_1^1(N^1)=S_1^1(m_2^1)\}} 1_{\{S_1^2(\bar{N}^2)=S_1^2(m_2^2)\}}$$
$$\times e^{-\Phi_\beta^2(\bar{N}^2,N^2)} 1_{\{S_1^2(\bar{N}^2)<S_1^2(\nu) \text{ for } \bar{N}^2<\nu\leq N^2\}}\big]$$

$$= \sum_{0<\bar{N}^2\leq N^2} \sum_{0\leq m_1^1<m_2^1\leq N^1} \sum_{0\leq m_1^2<m_2^2\leq \bar{N}^2} \sum_{y_0^1,y_0^2\in\mathbb{Z}^{d-1}} P_h[S^1(m_1^1)=(0,y_0^1-y^1)]$$

$$\times P_h[S^2(m_1^2)=(0,y_0^2-y^2)] \widetilde{B}_{h,\beta,y_0^1,y_0^2}^2(m_2-m_1)$$

$$\times P_h[S_1^1(N^1-m_2^1)=0] P_h[S_1^2(\bar{N}^2-m_2^2)=0] G_{h,\beta}^+(N^2-\bar{N}^2)$$

$$\leq G_{h,\beta}(N^1) \sum_{0 \leq m_1^1 < m_2^1 \leq N^1} \frac{P_h[S_1^1(m_1^1) = 0]}{G_{h,\beta}(m_1^1)} \frac{P_h[S_1^1(N^1 - m_2^1) = 0]}{G_{h,\beta}(N^1 - m_2^1)}$$

$$\times G_{h,\beta}(N^2) \sum_{0 \leq m_1^2 < \bar{m}_2^2 \leq N^2} \frac{P_h[S_1^2(m_1^2) = 0]}{G_{h,\beta}(m_1^2)} \frac{P_h[S_1^2(N^2 - \bar{m}_2^2) = 0]}{G_{h,\beta}(N^2 - \bar{m}_2^2)}$$

$$\times \sup_{\substack{n^1, M^2 \in \mathbb{N}, \\ y_0^1, y_0^2 \in \mathbb{Z}^{d-1}}} K_{h,\beta_0}^2 \sum_{n^2=1}^{M^2} \widetilde{B}_{h,\beta,y_0^1,y_0^2}^2(n) G_{h,\beta}^+(M^2 - n^2) e^{m_G(h,\beta)(n^1+M^2)},$$

where $\bar{m}_2^2$ replaces $m_2^2 + N^2 - \bar{N}^2$, $M^2$ stands for $\bar{m}_2^2 + m_1^2$, $n^1 \stackrel{\text{def}}{=} m_2^1 - m_1^1$, $n^2$ replaces $m_2^2 + m_1^2$ and $n \stackrel{\text{def}}{=} (n^1, n^2)$, and where (3.1) and Corollary 3.6 were used to split, respectively bound, $G_{h,\beta}(N^1)$ and $G_{h,\beta}(N^2)$.

It remains to bound the supremum in the above formula. By Lemma 4.2,

$$\widetilde{B}_{h,\beta,y_0^1,y_0^2}^2(n) e^{m_G(h,\beta)(n^1+n^2)}$$

$$\leq \sum_{k \in \mathbb{N}_0} \sum_{m \in \mathbb{N}_0} E_{h,\beta,y_0^1,y_0^2}[e^{\varphi_\beta(1) \sum_{i=0}^k \rho_{\sigma_i}} 1_{\{\sigma_k \leq m < \sigma_{k+1}\}} 1_{\{T(m)=n\}}]$$

$$\leq \sum_{k \in \mathbb{N}_0} \sum_{\tilde{n} \in \mathbb{N}_0^2: \tilde{n} \leq n} \sum_{m \in \mathbb{N}_0} \sum_{0 \leq \tilde{m} \leq m}$$

$$E_{h,\beta,y_0^1,y_0^2}[e^{\varphi_\beta(1) \sum_{i=0}^k \rho_{\sigma_i}} 1_{\{\sigma_k = \tilde{m}\}} 1_{\{T(\tilde{m}) = \tilde{n}\}}] P_{h,\beta}[T(m - \tilde{m}) = n - \tilde{n}]$$

$$= \sum_{k \in \mathbb{N}_0} \sum_{\tilde{n} \in \mathbb{N}_0^2: \tilde{n} \leq n} \sum_{\tilde{m} \in \mathbb{N}_0}$$

$$E_{h,\beta,y_0^1,y_0^2}[e^{\varphi_\beta(1) \sum_{i=0}^k \rho_{\sigma_i}} 1_{\{\sigma_k = \tilde{m}\}} 1_{\{T(\tilde{m}) = \tilde{n}\}}] \sum_{\hat{m} \in \mathbb{N}_0} P_{h,\beta}[T(\hat{m}) = n - \tilde{n}]$$

$$= \sum_{k \in \mathbb{N}_0} \sum_{\tilde{n} \in \mathbb{N}_0^2: \tilde{n} \leq n} E_{h,\beta,y_0^1,y_0^2}[e^{\varphi_\beta(1) \sum_{i=0}^k \rho_{\sigma_i}} 1_{\{\sigma_k < \infty\}} 1_{\{T(\sigma_k) = \tilde{n}\}}]$$

$$\times B_{h,\beta}^2(n - \tilde{n}) e^{m_G(h,\beta)(n^1 - \tilde{n}^1 + n^2 - \tilde{n}^2)},$$

where $B_{h,\beta}^2(\hat{n}) \stackrel{\text{def}}{=} E_h[\exp(-\Phi_\beta^{(2)}(\hat{n})); \mathrm{Br}^2(\hat{n})]$ for $\hat{n} \in \mathbb{N}_0^2$. We thus have

$$\sum_{n^2=1}^{M^2} \widetilde{B}_{h,\beta,y_0^1,y_0^2}^2(n) G_{h,\beta}^+(M^2 - n^2) e^{m_G(h,\beta)(n^1+M^2)}$$

$$\leq \sum_{k \in \mathbb{N}_0} \sum_{\tilde{n} \in \mathbb{N}_0^2: \tilde{n} \leq (n^1, M^2)} E_{h,\beta,y_0^1,y_0^2}[e^{\varphi_\beta(1) \sum_{i=0}^k \rho_{\sigma_i}} 1_{\{\sigma_k < \infty\}} 1_{\{T(\sigma_k) = \tilde{n}\}}]$$

$$\times \sup_{\hat{M}^2 \in \mathbb{N}_0} \sum_{\hat{n}^2=1}^{\hat{M}^2} B_{h,\beta}^2(\hat{n}) G_{h,\beta}^+(\hat{M}^2 - \hat{n}^2) e^{m_G(h,\beta)(\hat{n}^1 + \hat{M}^2)},$$

where $M^2$ stands for $M^2 - \tilde{n}^2$, $\hat{n}^1 \stackrel{\text{def}}{=} n^1 - \tilde{n}^1$, $\hat{n}^2$ stands for $n^2 - \tilde{n}^2$ and $\hat{n} \stackrel{\text{def}}{=} (\hat{n}^1, \hat{n}^2)$. By the independence of $\mathcal{S}^1$ and $\mathcal{S}^1$, we moreover have

$$B_{h,\beta}^2(\hat{n}) G_{h,\beta}^+(\hat{M}^2 - \hat{n}^2)$$

$$= \sum_{L=1}^{\infty} E_h[e^{-\Phi_\beta^1(\hat{n}^1) - \Phi_\beta^2(\hat{n}^2) - \Phi_\beta^2(\hat{n}^2, \hat{M}^2)} 1_{\{\mathcal{S}^1[\hat{n}^1] \text{ is a bridge of span } L\}}$$

$$\times 1_{\{\mathcal{S}^2[\hat{n}^2] \text{ is a bridge of span } L\}} 1_{\{L < S_1^2(\mu) \text{ for } \hat{n}^2 < \nu \leq \hat{M}^2\}}]$$

$$= \sum_{L=1}^{\infty} E_h[e^{-\Phi_\beta^1(\hat{n}^1)} 1_{\{\mathcal{S}^1[\hat{n}^1] \text{ is a bridge of span } L\}}]$$

$$\times E_h[e^{-\Phi_\beta^2(\hat{M}^2)} 1_{\{0 < S_1^2(\nu) \leq S_1^2(\hat{n}^2) = L < S_1^2(\nu) \text{ for } 0 < \mu \leq \hat{n}^2 < \nu \leq \hat{M}^2\}}].$$

Therefore, and by Proposition 3.2, we obtain

$$\sum_{\hat{n}^2=1}^{\hat{M}^2} B_{h,\beta}^2(\hat{n}) G_{h,\beta}^+(\hat{M}^2 - \hat{n}^2) \leq B_{h,\beta}(\hat{n}^1) G_{h,\beta}^+(\hat{M}^2) \leq e^{-m_G(h,\beta)(\hat{n}^1 + \hat{M}^2)},$$

which completes the proof of the Proposition 4.3. □

The next lemma gives a bound for the decay rate of the probability for large values of $\sigma_1$. It goes back to an estimate for the concentration of sums of independent, identically distributed random vectors with values in $\mathbb{Z}^{d-1}$ and to the fact that $\rho_{m+1}$ is zero when

$$\tau_{m+1}^1 + \tau_{m+1}^2 \leq \|\eta_m^1 - \eta_m^2\|_1.$$

Here, the dimension $d$ of the lattice comes into play explicitly for the first time.

LEMMA 4.4. *Suppose that $h > 0$ and $\beta_0 < \beta_c(h)$. Then, for any $\varepsilon > 0$, there exists $m_\varepsilon \in \mathbb{N}$ such that*

$$\sup_{y^1, y^2 \in \mathbb{Z}^{d-1}} P_{h,\beta,y^1,y^2}[\sigma_1 = m+1] \leq m^{-(d-1)/2+\varepsilon}$$

*for all $m \geq m_\varepsilon$ and $\beta \leq \beta_0$.*

PROOF. We distinguish whether the distance between the two random walks at the end of the $m$th strip is smaller or lager than $a_m \stackrel{\text{def}}{=} (\log m)^2$. For every pair of starting sites $y^1, y^2 \in \mathbb{Z}^{d-1}$, we have

$$P_{h,\beta,y^1,y^2}[\sigma_1 = m+1]$$
$$\leq P_{h,\beta,y^1,y^2}[\|\eta_m^1 - \eta_m^2\|_1 \leq a_m] + \sup_{\|\bar{y}^1 - \bar{y}^2\|_1 > a_m} P_{h,\beta,\bar{y}^1,\bar{y}^2}[\zeta_1 > 0].$$

Now, with $\gamma$ being chosen according to Theorem 3.7, we have

(4.2) $$K_1 \stackrel{\text{def}}{=} \sup_{\beta \leq \beta_0} E_{h,\beta}[e^{(\gamma/2)(\tau_1^1 + \tau_1^2)}] < \infty,$$

by continuity going back to the exponential gap for irreducible bridges. For any $\bar{y}^1, \bar{y}^2 \in \mathbb{Z}^{d-1}$ with $\|\bar{y}^1 - \bar{y}^2\|_1 > a_m$, the exponential Markov inequality thus implies that

$$P_{h,\beta,\bar{y}^1,\bar{y}^2}[\zeta_1 > 0] \leq P_{h,\beta}[\tau_1^1 + \tau_1^2 \geq \|\bar{y}^1 - \bar{y}^2\|_1]$$
$$\leq E_{h,\beta}[e^{(\gamma/2)(\tau_1^1 + \tau_1^2)}] e^{(-\gamma/2)a_m}$$
$$\leq K_1 m^{(-\gamma/2)\log m}$$

for all $m \in \mathbb{N}$ and $\beta \leq \beta_0$.

It remains to find an estimate for

$$P_{h,\beta,y^1,y^2}[\|\eta_m^1 - \eta_m^2\|_1 \leq a_m] \leq P_{h,\beta}[\|\eta_m^1 - \eta_m^2 - (y^2 - y^1)\|_2 \leq a_m].$$

Since the random sequence $(\eta_i^1 - \eta_i^2)_{i \in \mathbb{N}}$ is a Markov process with independent increments, the corollary to Theorem 6.2 in [6] yields

$$\sup_{y^1,y^2 \in \mathbb{Z}^{d-1}} P_{h,\beta}[\|\eta_m^1 - \eta_m^2 - (y^2 - y^1)\|_2 \leq a_m]$$
$$\leq K_2 \left(\frac{m}{a_m^2}\right)^{-(d-1)/2} \chi_1(a_m)^{-(d-1)/2},$$

where $K_2$ is a constant depending only on the dimension $d-1$ and where

$$\chi_1(u) \stackrel{\text{def}}{=} \inf_{t \in \mathbb{R}^{d-1}: \|t\|_2 = 1} E_{h,\beta}[\langle X - Y, t\rangle^2 \mathbf{1}_{\{\|X-Y\|_2 \leq u\}}],$$

$X, Y$ being independent copies of $\eta_1^1 - \eta_1^2$. A direct calculation and symmetry properties imply that

$$\inf_{t \in \mathbb{R}^{d-1}: \|t\|_2 = 1} E_{h,\beta}[\langle X - Y, t\rangle^2] = 2(d-1) E_{h,\beta}[(\eta_1^1 - \eta_1^2)_1^2],$$

where $(\eta_1^1 - \eta_1^2)_1$ denotes the first component of the vector $\eta_1^1 - \eta_1^2$. By the Cauchy–Schwarz inequality for sums and the exponential Markov inequality applied to $\|X - Y\|_2$, we also have

$$E_{h,\beta}[\langle X - Y, t\rangle^2 1_{\{\|X-Y\|_2 > u\}}] \leq \|t\|_2^2 e^{-u\gamma/4} E_{h,\beta}[\|X-Y\|_2^2 e^{(\gamma/4)\|X-Y\|_2}]$$
$$\leq \|t\|_2^2 e^{-u\gamma/4} \frac{32}{\gamma^2} E_{h,\beta}[e^{(\gamma/2)\|X-Y\|_2}]$$
$$\leq \|t\|_2^2 e^{-u\gamma/4} \frac{32}{\gamma^2} E_{h,\beta}[e^{(\gamma/2)\|X\|_1}] E_{h,\beta}[e^{(\gamma/2)\|Y\|_1}].$$

By continuity coming from the exponential gap in Theorem 3.7, there consequently exists a further constant $K_3 < \infty$ such that

$$\sup_{y^1,y^2 \in \mathbb{Z}^{d-1}} P_{h,\beta,y^1,y^2}[\|\eta_m^1 - \eta_m^2\|_1 \leq a_m] \leq K_3\left(\frac{m}{a_m^2}\right)^{-(d-1)/2}$$

for all $m \in \mathbb{N}$ and $\beta \leq \beta_0$, which completes the proof. $\square$

PROOF OF THEOREM 4.1. By renewing the Markov chain at times $\rho_{\sigma_1}, \ldots, \rho_{\sigma_k}$, we obtain

$$E_{h,\beta,y_0^1,y_0^2}[e^{\varphi_\beta(1)\sum_{i=1}^k \rho_{\sigma_i}} 1_{\{\sigma_k < \infty\}}] \leq \sup_{y^1,y^2 \in \mathbb{Z}^{d-1}} E_{h,\beta,y^1,y^2}[e^{\varphi_\beta(1)\rho_{\sigma_1}} 1_{\{\sigma_1 < \infty\}}]^k.$$

Theorem 4.1 thus follows from Proposition 4.3 once we show that

(4.3) $$\sup_{y^1,y^2 \in \mathbb{Z}^{d-1}} \sup_{\beta \leq \beta_0} E_{h,\beta,y^1,y^2}[e^{\varphi_\beta(1)\rho_{\sigma_1}} 1_{\{\sigma_1 < \infty\}}] < 1$$

for $d \geq 4$ and $\beta_0 > 0$ small enough.

To this end, we choose $0 < \varepsilon < \frac{1}{2}$ and $p, q > 1$ with $\frac{1}{q}(\frac{3}{2} - \varepsilon) > 1$ and $\frac{1}{p} + \frac{1}{q} = 1$. By Lemma 4.4, there exists $m_\varepsilon \in \mathbb{N}$ such that the Hölder inequality implies

(4.4) 
$$E_{h,\beta,y^1,y^2}[e^{\varphi_\beta(1)\rho_{m+1}} 1_{\{\sigma_1 = m+1\}}]$$
$$\leq E_{h,\beta,y^1,y^2}[e^{p\varphi_\beta(1)\rho_{m+1}}]^{1/p} P_{h,\beta,y^1,y^2}[\sigma_1 = m+1]^{1/q}$$
$$\leq E_{h,\beta}[e^{p\varphi_\beta(1)(\tau_1^1+\tau_1^2+2)}]^{1/p} m^{(-1/q)((d-1)/2-\varepsilon)}$$
$$\leq (K_1 + e^{2\varphi_{\beta_0}(1)})^{1/p} m^{(-1/q)((d-1)/2-\varepsilon)}$$

for all $m \geq m_\varepsilon$ and $\beta \leq \beta_0$, where $K_1$ is defined in (4.2) and $\beta_0$ needs to be chosen small enough. The exponential Hölder inequality moreover yields

$$P_{h,\beta,y^1,y^2}[\sigma_1 = m+1] \leq P_{h,\beta,y^1,y^2}[\zeta_{m+1} > 0]$$
$$\leq P_{h,\beta}[T^1(m+1) + T^2(m+1) \geq \|y^1 - y^2\|_1]$$
$$\leq K_1^{m+1} e^{(-\gamma/2)\|y^1 - y^2\|_1}$$

for all $m \in \mathbb{N}$ and $\beta \leq \beta_0$, where $\gamma$ is chosen according to (4.2). Since $\frac{1}{q}(\frac{d}{2} - \varepsilon) > 1$ for $d \geq 4$, we thus obtain

$$(4.5) \qquad \lim_{k \to \infty} \sup_{y^1, y^2 : \|y^1 - y^2\|_1 > k} \sup_{\beta \leq \beta_0} E_{h,\beta,y^1,y^2}[e^{\varphi_\beta(1)\rho_{\sigma_1}} 1_{\{\sigma_1 < \infty\}}] = 0$$

for $\beta_0$ small enough.

Now, for any $k \in \mathbb{N}$ and $y^1 \neq y^2$, we have

$$P_{h,0,y^1,y^2}[\sigma_1 = \infty] \geq P_{h,0,y^1,y^2}[\zeta_1 = 0, \|\eta_1^1 - \eta_1^2\|_1 \geq k]$$
$$\times \inf_{\|\bar{y}^1 - \bar{y}^2\|_1 \geq k} P_{h,0,\bar{y}^1,\bar{y}^2}[\sigma_1 = \infty]$$

and similarly

$$P_{h,0}[\sigma_1 = \infty] \geq P_{h,0}[\zeta_1 = 2, \|\eta_1^1 - \eta_1^2\|_1 \geq k] \inf_{\|\bar{y}^1 - \bar{y}^2\|_1 \geq k} P_{h,0,\bar{y}^1,\bar{y}^2}[\sigma_1 = \infty].$$

Therefore, and by (4.5), we obtain

$$(4.6) \qquad \sup_{y^1, y^2 \in \mathbb{Z}^{d-1}} P_{h,0,y^1,y^2}[\sigma_1 < \infty] < 1.$$

Finally, for any $y^1, y^2 \in \mathbb{Z}^{d-1}$, observe that

$$E_{h,\beta,y^1,y^2}[e^{\varphi_\beta(1)\rho_{\sigma_1}} 1_{\{\sigma_1 < \infty\}}]$$

is continuous in $\beta = 0$, the continuity going back to (4.4) and the exponential gap in Theorem 3.7. Consequently, (4.5) and (4.6) now imply (4.3), completing the proof of the theorem. $\square$

**Acknowledgment.** I wish to thank Prof. Erwin Bolthausen for having given me the opportunity to study this fascinating topic within the framework of a doctoral thesis under his supervision. With his great mathematical intuition, he perceived a way to handle the model, and he gave me important support in the realization of my thesis.

## REFERENCES

[1] ALBEVERIO, S. and ZHOU, X. Y. (1996). A martingale approach to directed polymers in a random environment. *J. Theoret. Probab.* **9** 171–189. MR1371075
[2] BOLTHAUSEN, E. (1989). A note on the diffusion of directed polymers in a random environment. *Comm. Math. Phys.* **123** 529–534. MR1006293
[3] CARMONA, P. and HU, Y. (2002). On the partition function of a directed polymer in a Gaussian random environment. *Probab. Theory Related Fields* **124** 431–457. MR1939654
[4] CHAYES, J. T. and CHAYES, L. (1986). Ornstein–Zernike behavior for self-avoiding walks at all noncritical temperatures. *Comm. Math. Phys.* **105** 221–238. MR0849206

[5] COMETS, F., SHIGA, T. and YOSHIDA, N. (2003). Directed polymers in a random environment: Path localization and strong disorder. *Bernoulli* **9** 705–723. MR1996276

[6] ESSEN, C. G. (1968). On the concentration function of a sum of independent random variables. *Z. Wahrsch. Verw. Gebiete* **9** 290–308. MR0231419

[7] FLURY, M. (2007). Large deviations and phase transition for random walks in random nonnegative potentials. *Stochastic Process. Appl.* **117** 596–612. MR2320951

[8] GREVEN, G. and DEN HOLLANDER, F. (1992). Branching random walk in random environment: Phase transition for local and global growth rates. *Probab. Theory Related Fields* **91** 195–249. MR1147615

[9] IOFFE, D. (1998). Ornstein–Zernike behavior and analyticity of shapes for self-avoiding walks on $\mathbb{Z}^d$. *Markov Process. Related Fields* **4** 323–350. MR1670027

[10] IMBRIE, J. Z. and SPENCER, T. (1988). Diffusion of directed polymers in a random environment. *J. Statist. Phys.* **52** 609–626. MR0968950

[11] MADRAS, M. and SLADE, G. (1993). *The Self-Avoiding Walk.* Birkhäuser, Boston. MR1197356

[12] SZNITMAN, A. S. (1998). *Brownian Motion, Obstacles and Random Media.* Springer, Berlin. MR1717054

[13] TALAGRAND, M. (1996). A new look at independence. *Ann. Probab.* **24** 1–34. MR1387624

[14] TRACHSLER, M. (1999). Phase transitions and fluctuations for random walks with drift in random potentials. Ph.D. thesis, Univ. Zurich.

[15] ZERNER, M. P. W. (1998). Directional decay of the Green's function for a random nonnegative potential on $\mathbb{Z}^d$. *Ann. Appl. Probab.* **8** 246–280. MR1620370

MATHEMATISCHES INSTITUT
UNIVERSITÄT TÜBINGEN
AUF DER MORGENSTELLE 10
D-72076 TÜBINGEN
GERMANY
E-MAIL: mflury@amath.unizh.ch